\documentclass[english,reqno]{amsart}
\usepackage{times}
\usepackage[T1]{fontenc}
\usepackage{graphicx}
\usepackage{amssymb}

\makeatletter
 \theoremstyle{plain}    
 \newtheorem{thm}{Theorem}[section]
 \numberwithin{equation}{section} %% Comment out for sequentially-numbered
 \numberwithin{figure}{section} %% Comment out for sequentially-numbered
 \theoremstyle{plain}
 \theoremstyle{definition}
  \newtheorem*{example*}{Example}
 \theoremstyle{plain}    
 \newtheorem{prop}[thm]{Proposition} %%Delete [thm] to re-start numbering
 \theoremstyle{plain}    
 \newtheorem{cor}[thm]{Corollary} %%Delete [thm] to re-start numbering
 \theoremstyle{remark}
 \newtheorem*{rem*}{Remark}
 \theoremstyle{plain}    
 \newtheorem{lem}[thm]{Lemma} %%Delete [thm] to re-start numbering

\usepackage{psfrag}
 \usepackage{bbm}

\newcommand{\Q}{\mathbb {Q}}

\newcommand{\I}{\mathbb {I}}

\AtBeginDocument{
  
}

\usepackage{babel}
\makeatother
\begin{document}
\title{Critical waiting time processes in infinite ergodic theory}

\author{Marc Kesseböhmer and Mehdi Slassi}

\address{Fachbereich 3 - Mathematik und Informatik, Universit\"{a}t Bremen,
D--28359 Bremen, Germany}

\email{mhk@math.uni-bremen.de, slassi@math.uni-bremen.de}

\begin{abstract}
We study limit laws for return time processes defined on infinite
conservative ergodic measure preserving dynamical systems. Especially
for the critical cases with purely atomic limiting distribution we
derive distorted processes posessing non-degenerated limits. For these
processes also large deviation asymptotics are stated. The Farey map
is used as an illustrating example giving new insides into the metric
number theory of continued fractions.
\end{abstract}
\maketitle

\section{Introduction and preliminaries}

\let\languagename\relax

Throughout this paper $\left(X,T,\mathcal{A},\mu\right)$ will denote
a conservative ergodic measure preserving dynamical systems where
$\mu$ is an infinite $\sigma$-finite measure. Kac's Theorem implies
that in this situation the mean return time to sets of finite positive
measure is infinite. In terms of Markov chains this corresponds to
the \emph{null recurrent} setting. Hence, new probabilistic properties
of such dynamical systems lead to interesting results for null recurrent
Markov chains, whereas known results for these Markov chains sometimes
allow analog statements for infinite measure preserving transformations. 

In this paper we are going to study the critical cases of the generalized
Thaler's Dynkin-Lamperti laws describing the asymptotic behaviour
of the following processes.

\begin{itemize}
\item $Z_{n}(x):=\left\{ \begin{array}{ll}
\max\{ k\leq n:\; T^{k}(x)\in A\}, & x\in A_{n}:=\bigcup_{k=0}^{n}T^{-k}A,\\
0, & \textrm{else.}\end{array}\right.$
\item $Y_{n}\left(x\right):=\min\left\{ k>n:\; T^{k}\left(x\right)\in A,\right\} \quad x\in X,$
\item $V_{n}\left(x\right):=Y_{n}-Z_{n}.$
\end{itemize}
Namely, it is shown in \cite{Thaler:98} that $Z_{n}/n$, $Y_{n}/n$,
$V_{n}/N$ all converge strongly in distribution to certain random
variables depending only on the exponent $\left(1-\alpha\right)$
of the wandering rate (cf. (T), Subsection \ref{sub:Limit-laws.}).
For certain values of $\alpha$ these random variables turn out to
be atomic. In order to derive non-degenerated results also for these
cases we consider \emph{distorted processes, i.e.} \[
\frac{F\left(X_{n}\right)}{F\left(n\right)}\quad\textrm{and }\quad\frac{F\left(n-X_{n}\right)}{F\left(n\right)},\]
where $F$ is a regularly varying function and $\left(X_{n}\right)$
denotes any of the sequences $\left(Y_{n}\right)$, $\left(V_{n}\right)$.
In particular we introduce the processes\[
\Lambda_{n}:=\frac{\sum_{k=0}^{V_{n}}\mu\left(A\cap\left\{ \varphi>k\right\} \right)}{\mu\left(A_{n}\right)},\qquad\Gamma_{n}:=\frac{n\sum_{k=0}^{V_{n}}\mu\left(A\cap\left\{ \varphi>k\right\} \right)}{\mu\left(A_{n}\right)V_{n}},\]
\[
\Delta_{n}:=\frac{\sum_{k=0}^{Y_{n}-n}\mu\left(A\cap\left\{ \varphi>k\right\} \right)}{\mu\left(A_{n}\right)},\qquad\Theta_{n}:=\frac{n\sum_{k=0}^{Y_{n}}\mu\left(A\cap\left\{ \varphi>k\right\} \right)}{\mu\left(A_{n}\right)Y_{n}}.\]
 We call $\Lambda_{n}$ the \emph{distorted total waiting time process}
and $\Delta_{n}$ the \emph{distorted residual waiting time process}.
In here, \begin{equation}
\varphi(x)=\inf\{ n\geq1:\; T^{n}(x)\in A\},\;\; x\in X,\label{returntime}\end{equation}
 denotes the \emph{first return time} to the set $A$.

We remark that the analog questions for $\left(Z_{n}\right)$ are
already treated in \cite{KesseboehmerSlassi:05}. In \cite{KessboehmerSlassi:05}
some limit theorems for $\left(Z_{n}\right)$ have been applied to
the Farey interval map deriving new number theoretical results for
continued fractions. In the last section of this paper we also develop
some consequences of the main theorems for continued fractions. 

Finally, we would like to point out that other related results can
be found in \cite{ThalerZweimueller:06}, \cite{Zweimueller:03}.

\subsection{Infinite ergodic theory. \label{sub:Infinite-ergodic-theory} }

A characterization of $\left(X,T,\mathcal{A},\mu\right)$ being a
con\-ser\-vative ergodic measure preserving dynamical system where
$\mu$ is an infinite $\sigma$-finite measure as used in this paper
will be given in terms of the transfer operator below. For further
definitions and details we refer the reader to \cite{Aaronson:97}.

Let \[
\mathcal{P}_{\mu}:=\left\{ \nu:\nu\:\textrm{probability measure on}\:\mathcal{A}\,\textrm{with }\nu\ll\mu\right\} \]
 denote the set of probability measures on $\mathcal{A}$ which are
absolutely continuous with respect to $\mu$. The measures from $\mathcal{P}_{\mu}$
represent the admissible initial distributions for the processes under
consideration. With $\mathcal{P}_{\mu}$ we will sometimes also denote
the set of the corresponding densities. 

Let us recall the notion of the wandering rate. For a fixed set $A\in\mathcal{A}$
with $0<\mu\left(A\right)<\infty$ we set 

\[
A_{n}:=\bigcup_{k=0}^{n}T^{-k}A\quad\mathrm{and}\quad W_{n}:=W_{n}\left(A\right):=\mu\left(A_{n}\right),\qquad n\geq0,\]
and call the sequence $\left(W_{n}\left(A\right)\right)$ the \emph{wandering
rate} of $A.$ The following identities hold\[
W_{n}\left(A\right)=\sum_{k=0}^{n}\mu\left(A\cap\{\varphi>k\}\right)=\int_{A}\min(\varphi,n+1)\, d\mu,\qquad n\geq0.\]
 Since $T$ is conservative and ergodic, for all $\nu\in\mathcal{P_{\mu}}$,
we have\[
\nu\left(\left\{ \varphi<\infty\right\} \right)=1,\;\nu\left(\left\{ Y_{n}<\infty\right\} \right)=1\;\textrm{for\; all}\; n\geq1,\,\quad\textrm{and}\quad\lim_{n\to\infty}\nu\left(A_{n}\right)=1\]

To explore the stochastic properties of a non-singular transformation
of a $\sigma$--finite measure space it is often useful to study the
long-term behaviour of the iterates of its \emph{transfer operator}
\[
\hat{T}:L_{1}\left(\mu\right)\longrightarrow L_{1}\left(\mu\right),\; f\longmapsto\hat{T}\left(f\right):=\frac{d\left(\nu_{f}\circ T^{-1}\right)}{d\mu},\]
where $\nu_{f}$ denote the measure with density $f$ with respect
to $\mu$. Clearly, $\hat{T}$ is a positive linear operator characterized
by\[
\int_{B}\hat{T}\left(f\right)\; d\mu=\int_{T^{-1}\left(B\right)}f\; d\mu,\qquad f\in L_{1}\left(\mu\right),\quad B\in\mathcal{A}.\]
An approximation argument shows that equivalently for all $f\in L_{1}\left(\mu\right)$
and $g\in L_{\infty}\left(\mu\right)$\[
\int_{X}\hat{T}\left(f\right)\cdot g\; d\mu=\int_{X}f\cdot g\circ T\; d\mu.\]
The ergodic properties of $(X,T,\mathcal{A},\mu)$ can be characterized
in terms of the transfer operator in the following way (cf. \cite{Aaronson:97}).
A system is conservative and ergodic if and only if for all $f\in L_{1}^{+}\left(\mu\right):=\left\{ f\in L_{1}\left(\mu\right):\; f\geq0\;\mathrm{and}\;\int_{X}f\; d\mu>0\right\} $
we have $\mu$-a.e. \[
\sum_{n\geq0}\hat{T}^{n}\left(f\right)=\infty.\]
 Note that invariance of $\mu$ under $T$ means $\hat{T}\left(1\right)=1.$

The following two definitions are in many situation crucial within
infinite ergodic theory.

\begin{itemize}
\item A set $A\in\mathcal{A}$ with $0<\mu\left(A\right)<\infty$ is called
\emph{uniform for} $f\in\mathcal{P_{\mu}}$ if there exists a sequence
$\left(b_{n}\right)$ of positive reals such that\[
\frac{1}{d_{n}}\sum_{k=0}^{n-1}\hat{T}^{k}\left(f\right)\;\longrightarrow\;1\qquad\mu-\textrm{a.e. \; uniformly\; on}\; A\]
(i.e. uniform convergence in $L_{\infty}\left(\mu|_{A\cap\mathcal{A}}\right)$).
\item The set $A$ is called a \emph{uniform} set \emph{}if it is uniform
for some $f\in\mathcal{P_{\mu}}.$
\end{itemize}
Note that from \cite{Aaronson:97}, Proposition 3.8.7, we know, that
$\left(b_{n}\right)$ is regularly varying with exponent $\alpha$
(for the definition of this property see Subsection \ref{sub:Classical-results-on})
if and only if $\left(W_{n}\right)$ is regularly varying with exponent
$\left(1-\alpha\right)$. In this case $\alpha$ lies in the interval
$\left[0,1\right]$ and \begin{equation}
d_{n}W_{n}\sim\frac{n}{\Gamma\left(1+\alpha\right)\Gamma\left(2-\alpha\right)}.\label{eq:AsymptAaronson}\end{equation}
In here, $c_{n}\sim b_{n}$ for some sequences $\left(c_{n}\right)$
and $\left(b_{n}\right)$ means that $b_{n}\not=0$ has only finitely
many exceptions and $\lim_{n\to\infty}\frac{c_{n}}{b_{n}}=1$.

Next, we recall the notion of uniformly returning sets, which will
be used to state the conditions in Theorem \ref{theo2} (cf. \cite{KesseboehmerSlassi:05},
Subsection1.2)

\begin{itemize}
\item A set $A\in\mathcal{A}$ with $0<\mu\left(A\right)<\infty$ is called
\emph{uniformly returning} \emph{for} $f\in\mathcal{P_{\mu}}$ if
there exists an positive increasing sequence $\left(b_{n}\right)$
such that \[
b_{n}\hat{T}^{n}\left(f\right)\;\longrightarrow\;1\quad\mu-\textrm{a.e. \; uniformly\; on}\; A.\]

\item The set $A$ is \emph{}called \emph{uniformly} \emph{returning} if
it is uniformly returning for some $f\in\mathcal{P_{\mu}}.$\medbreak
\end{itemize}
From \cite{KesseboehmerSlassi:05}, Proposition 1.2, we know that
for $\beta\in\left[0,1\right)$ we have that $\left(b_{n}\right)$
is regularly varying with exponent $\beta$ if and only if $\left(W_{n}\right)$
is regularly varying with the same exponent. In this case,\[
b_{n}\sim W_{n}\Gamma\left(1-\beta\right)\Gamma\left(1+\beta\right)\qquad\left(n\to\infty\right).\]
 Also (see \cite{KesseboehmerSlassi:05}, Proposition 1.1) we have
that every uniformly returning set is uniform. In the proof of this
fact it is shown that there exists a funcition $f\in\mathcal{P}_{\mu}$
such that $A$ is uniform as well as uniformly returning for $f$.
This observation will be relevant in Theorem \ref{thm:neue2}. Under
some extra conditions the reverse implication is also true (cf. \cite{KessboehmerSlassi:05},
Proposition 2.6). 

\begin{example*}
Let $T:\left[0,1\right]\longrightarrow\left[0,1\right]$ be an interval
map with two increasing full branches and an indifferent fixed point
at $0$ satisfying Thaler's conditions in \cite{Thaler:00}. Then
any set $A\in\mathcal{B}_{\left[0,1\right]}$ with positive distance
from the indifferent fixed point $0$ and positive Lebesgue measure
$\lambda\left(A\right)$ is uniformly returning.
\end{example*}
Sometimes the limiting behaviour of processes defined in terms of
a non-singular transformation does not depend on the initial distribution.
This is formalized as follows.

Let $\nu$ be a probability measure on the measurable space $(X,\mathcal{A})$
and $\left(R_{n}\right)_{n\geq1}$ be a sequence of measurable real
functions on $X$. Then distributional convergence of $\left(R_{n}\right)_{n\geq1}$
w.r.t. $\nu$ to some random variable $R$ with values in $\left[-\infty,\infty\right]$
will be denoted by $R_{n}\stackrel{\nu}{\Longrightarrow}R$. \emph{Strong
distributional convergence} abbreviated by $R_{n}\stackrel{\mathcal{L\left(\mu\right)}}{\Longrightarrow}R$
on the $\sigma$--finite measures space $\left(X,\mathcal{A,\mu}\right)$
means that $R_{n}\stackrel{\nu}{\Longrightarrow}R$ for all $\nu\in\mathcal{P_{\mu}}$.
In particular for $c\in\left[-\infty,\infty\right]$,\[
R_{n}\stackrel{\mathcal{L\left(\mu\right)}}{\Longrightarrow}c\quad\Longleftrightarrow\quad R_{n}\longrightarrow c\textrm{\quad locally in measure,}\]
which we also denote by $R_{n}\stackrel{\mu}{\longrightarrow}c$.

\section{Statements of main results.\label{sub:Limit-laws.}}

We begin this section with recalling the following interesting limit
laws for the processes $Y_{n}$ and $V_{n}$ which are due to Thaler
\cite{Thaler:98}.

\begin{itemize}
\item [(T)] \textbf{Thaler's Dynkin-Lamperti Theorem.} Let $A\in\mathcal{A}$
with $0<\mu\left(A\right)<\infty$ be a uniform set. If the wandering
rate $\left(W_{n}\left(A\right)\right)$ is regularly varying with
exponent $1-\alpha$ for $\alpha\in\left[0,1\right]$, then we have

\begin{itemize}
\item [{\rm (1)}]\[
\frac{Y_{n}}{n}\;\stackrel{\mathcal{L\left(\mu\right)}}{\Longrightarrow}\;\varphi_{\alpha},\]
where $\varphi_{\alpha}$ , $\alpha\in\left(0,1\right)$ , denotes
the random variable on $\left[1,\infty\right)$ with density \[
f_{\varphi_{\alpha}}\left(x\right)=\frac{\sin\pi\alpha}{\pi}\frac{1}{x\left(x-1\right)^{\alpha}},\quad x>1,\]
and $\varphi_{0}=\infty$, $\varphi_{1}=1.$
\item [{\rm (2)}]\[
\frac{V_{n}}{n}\;\stackrel{\mathcal{L\left(\mu\right)}}{\Longrightarrow}\;\eta_{\alpha},\]
where $\eta_{\alpha}$ , $\alpha\in\left(0,1\right)$ , denotes the
random variable on $\left[0,\infty\right)$ with density \[
f_{\eta_{\alpha}}\left(x\right)=\frac{\sin\pi\alpha}{\pi}\frac{1-\left(\max\left\{ 1-x,0\right\} \right)^{\alpha}}{x^{1+\alpha}},\qquad x>0,\]
and $\eta_{0}=\infty,$ $\eta_{1}=0$.
\end{itemize}
\end{itemize}
To apply (T) to the distorted processes we need the following proposition
from \cite{KessboehmerSlassi:05}. Its first part was independently
proved by Thaler in \cite{Thaler:05}. 

\begin{prop}
\label{propo0} Let $\left(\Omega,\mathcal{F},\mathbb{P}\right)$
be a probability space, let $Y_{n}:\Omega\longrightarrow\left[0,\infty\right]$
be measurable $\left(n\geq1\right)$, and let $Y$ be a random variable
with values in $\left[0,\infty\right]$.
\begin{itemize}
\item [{\rm (1)}] If $\mathbb{P}\left(Y=0\right)=0=\mathbb{P}\left(Y=\infty\right)$
and $F$ is a regularly varying function with exponent $\beta\in\mathbb{R},$
then \[
\frac{Y_{n}}{n}\stackrel{\mathbb{P}}{\Longrightarrow}Y\qquad\quad\implies\qquad\quad\frac{F\left(Y_{n}\right)}{F\left(n\right)}\stackrel{\mathbb{P}}{\Longrightarrow}Y^{\beta}.\]

\item [{\rm (2)}] If $Y=0$ and $F$ is a regularly varying function with
exponent $\beta\in\mathbb{R}\setminus\left\{ 0\right\} $ then\[
\frac{Y_{n}}{n}\stackrel{\mathbb{P}}{\Longrightarrow}0\qquad\quad\implies\qquad\quad\frac{F\left(Y_{n}\right)}{F\left(n\right)}\stackrel{\mathbb{P}}{\Longrightarrow}\left\{ \begin{array}{ll}
0 & \textrm{for }\beta>0\\
\infty & \textrm{for }\beta<0\end{array}\right..\]

\end{itemize}
\end{prop}
\begin{itemize}
\item [{\rm (3)}] If $Y=\infty$ and $F$ is a regularly varying function
with exponent $\beta\in\mathbb{R}\setminus\left\{ 0\right\} $ then\[
\frac{Y_{n}}{n}\stackrel{\mathbb{P}}{\Longrightarrow}\infty\qquad\quad\implies\qquad\quad\frac{F\left(Y_{n}\right)}{F\left(n\right)}\stackrel{\mathbb{P}}{\Longrightarrow}\left\{ \begin{array}{ll}
\infty & \textrm{for }\beta>0\\
0 & \textrm{for }\beta<0\end{array}\right..\]

\end{itemize}
The following corollary is a direct consequence of (T), Proposition\ref{propo0},
and the fact that \[
\Lambda_{n}=\frac{F\left(V_{n}\right)}{F\left(n\right)},\quad,\Gamma_{n}:=\frac{G\left(V_{n}\right)}{G\left(n\right)},\quad\Delta_{n}=\frac{F\left(Y_{n}-n\right)}{F\left(n\right)},\quad\Theta_{n}=\frac{G\left(Y_{n}\right)}{G(n)},\]
 with $F\left(n\right):=W_{n}$ and $G\left(n\right):=W_{n}/n$. 

\begin{cor}
\label{coro 1} Let $A\in\mathcal{A}$ with $0<\mu\left(A\right)<\infty$
be a uniform set such that the wandering rate $\left(W_{n}\right)$
is regularly varying with exponent $1-\alpha$. 
\begin{enumerate}
\item If $0\leq\alpha<1,$ then we have\[
\Lambda_{n}\;\stackrel{\mathcal{L\left(\mu\right)}}{\Longrightarrow}\;\lambda_{\alpha},\]
where $\zeta_{\alpha}$ denotes the random variable on $\left[0,\infty\right)$
with density \[
f_{\zeta_{\alpha}}\left(x\right)=\frac{1}{1-\alpha}\frac{\sin\pi\alpha}{\pi}\frac{1-\left(\max\left\{ 1-x^{\frac{1}{1-\alpha}},0\right\} \right)^{\alpha}}{x^{\frac{1}{1-\alpha}}},\qquad\alpha\in\left(0,1\right),\]
and $\lambda_{0}=\infty$ (cf. Fig. \ref{cap:Figdensity1}).%
\begin{figure}
\psfrag{0.8}{{\small $\alpha=0.98$}}

\psfrag{0.6}{{\small $\alpha=0.5$}}

\psfrag{0.4}{{\small $\alpha=0.4$}}

\psfrag{0.2}{{\small $\alpha=0.2$}}\psfrag{0.22}{{\small $\alpha=0.005$}} \psfrag{0}{{\small $0$}}\psfrag{1}{{\small $1$}}\psfrag{0.5}{{\small $0.5$}}\psfrag{1.5}{{\small $1$}}\psfrag{1.5}{{\small $1.5$}}\psfrag{2}{{\small $2$}}\psfrag{2.5}{{\small $2.5$}}\psfrag{3}{{\small $3$}}

\includegraphics[%
  width=0.70\textwidth]{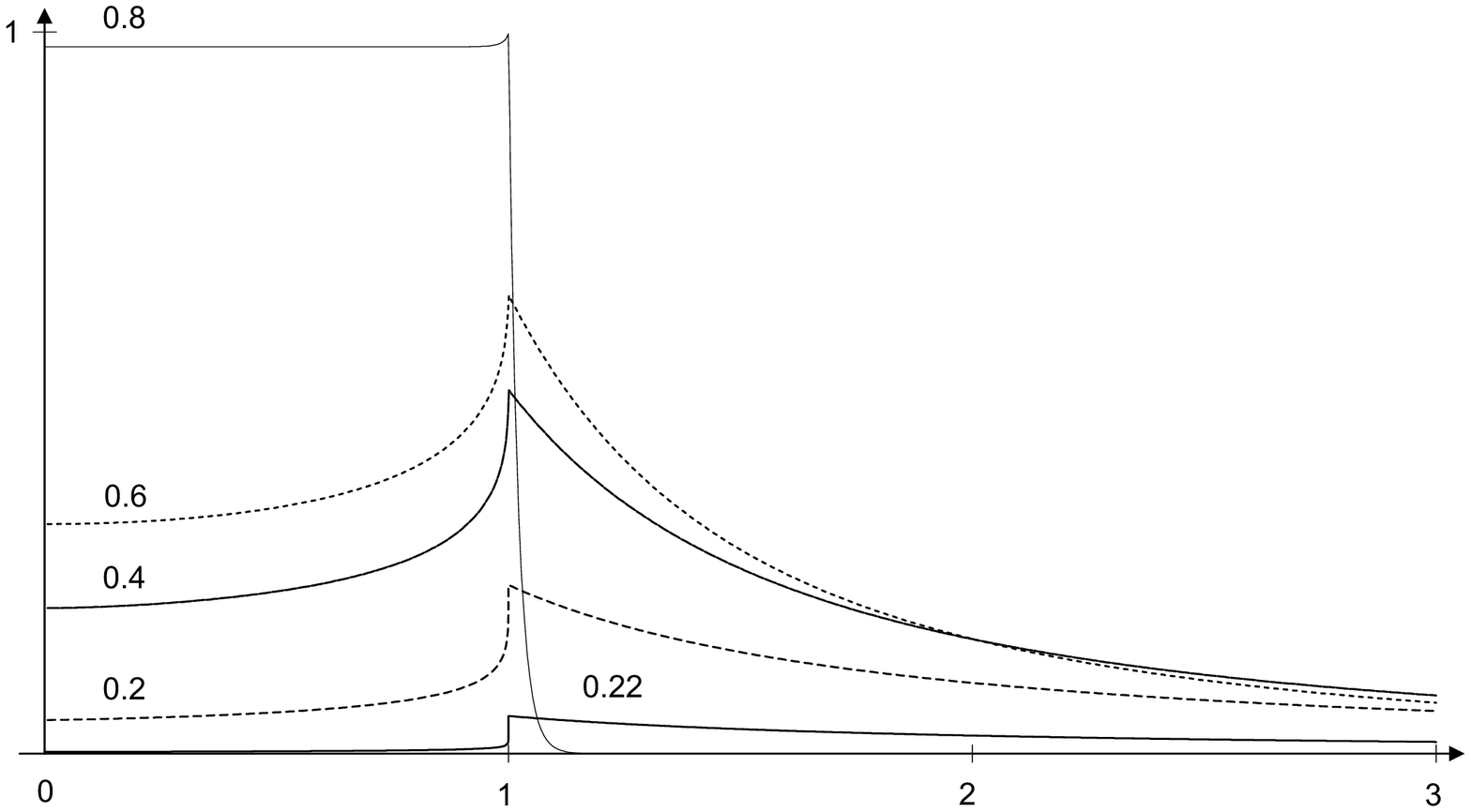}

\caption{\label{cap:Figdensity1} The densities $f_{\lambda_{\alpha}}$ of
the limiting distribution of the normalized Kac process for different
values of $\alpha\in\left(0,1\right)$. The extreme case is given
by $\lambda_{0}=\infty$. }
\end{figure}

\item If $0<\alpha\leq1,$ then we have\[
\Gamma_{n}\;\stackrel{\mathcal{L\left(\mu\right)}}{\Longrightarrow}\;\gamma_{\alpha},\]
where $\gamma_{\alpha}$ denotes the random variable on $\left[0,\infty\right)$
with density \[
f_{\gamma_{\alpha}}\left(x\right)=\frac{\sin\pi\alpha}{\alpha\pi}\left(1-\left(\max\left\{ 1-x^{\frac{-1}{\alpha}},0\right\} \right)^{\alpha}\right),\qquad\alpha\in\left(0,1\right),\]
and $\gamma_{1}=\infty$ (cf. Fig. \ref{cap:Figdensity4}).%
\begin{figure}
\psfrag{0.98}{{\small $\alpha=0.98$}}

\psfrag{0.5}{{\small $\alpha=0.5$}}\psfrag{0.7}{{\small $\alpha=0.7$}}

\psfrag{0.3}{{\small $\alpha=0.3$}}

\psfrag{0.2}{{\small $\alpha=0.2$}}\psfrag{0.05}{{\small $\alpha=0.05$}} \psfrag{0}{{\small $0$}}\psfrag{1}{{\small $1$}}\psfrag{4}{{\small $4$}}\psfrag{1.5}{{\small $1.5$}}\psfrag{2}{{\small $2$}}\psfrag{2.5}{{\small $2.5$}}\psfrag{3}{{\small $3$}}

\includegraphics[%
  width=0.70\textwidth]{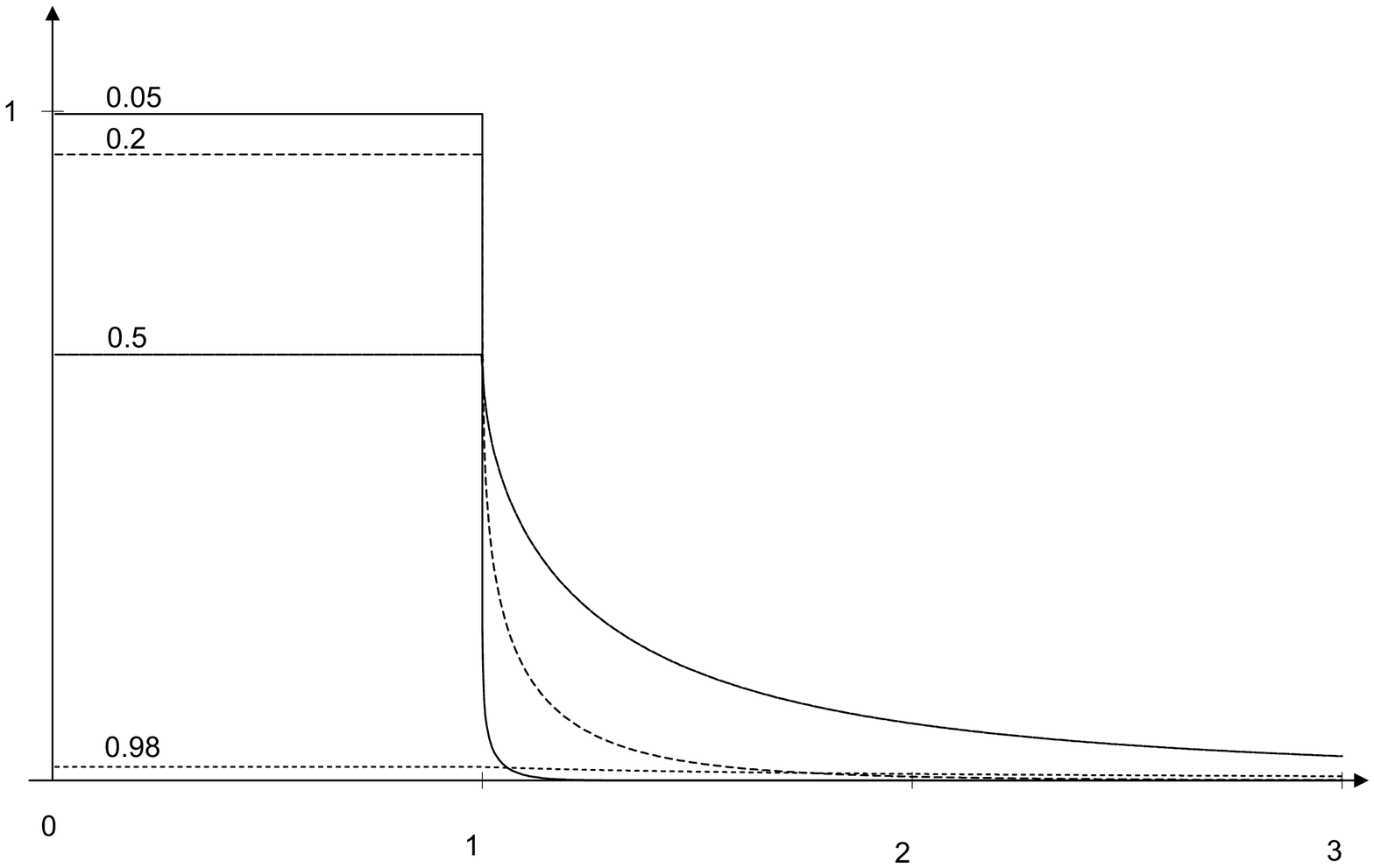}

\caption{\label{cap:Figdensity4} The densities $f_{\gamma_{\alpha}}$ of
the limiting distribution of the normalized Kac process for different
values of $\alpha\in\left(0,1\right)$. The extreme case is given
by $\gamma_{1}=\infty$. }
\end{figure}

\item If $0\leq\alpha<1,$ then we have\[
\Delta_{n}\;\stackrel{\mathcal{L\left(\mu\right)}}{\Longrightarrow}\;\delta_{\alpha},\]
where $\gamma_{\alpha}$ denotes the random variable on $\left[0,\infty\right)$
with density \[
f_{\delta_{\alpha}}\left(x\right)=\frac{1}{1-\alpha}\frac{\sin\pi\alpha}{\pi}\frac{1}{1+x^{\frac{1}{1-\alpha}}},\qquad\alpha\in\left(0,1\right),\]
and $\delta_{0}=\infty$. (cf. Fig. \ref{cap:Figdensity3}). %
\begin{figure}
\psfrag{0.99}{{\small $\alpha=0.99$}}\psfrag{0.7}{{\small $\alpha=0.7$}}

\psfrag{0.5}{{\small $\alpha=0.5$}}

\psfrag{0.3}{{\small $\alpha=0.3$}}

\psfrag{0.2}{{\small $\alpha=0.2$}}\psfrag{0.05}{{\small $\alpha=0.05$}} \psfrag{0}{{\small $0$}}\psfrag{1}{{\small $1$}}\psfrag{1.5}{{\small $1$}}\psfrag{1.5}{{\small $1.5$}}\psfrag{2}{{\small $2$}}\psfrag{2.5}{{\small $2.5$}}\psfrag{3}{{\small $3$}}\psfrag{4}{{\small $4$}}

\includegraphics[%
  width=0.70\textwidth]{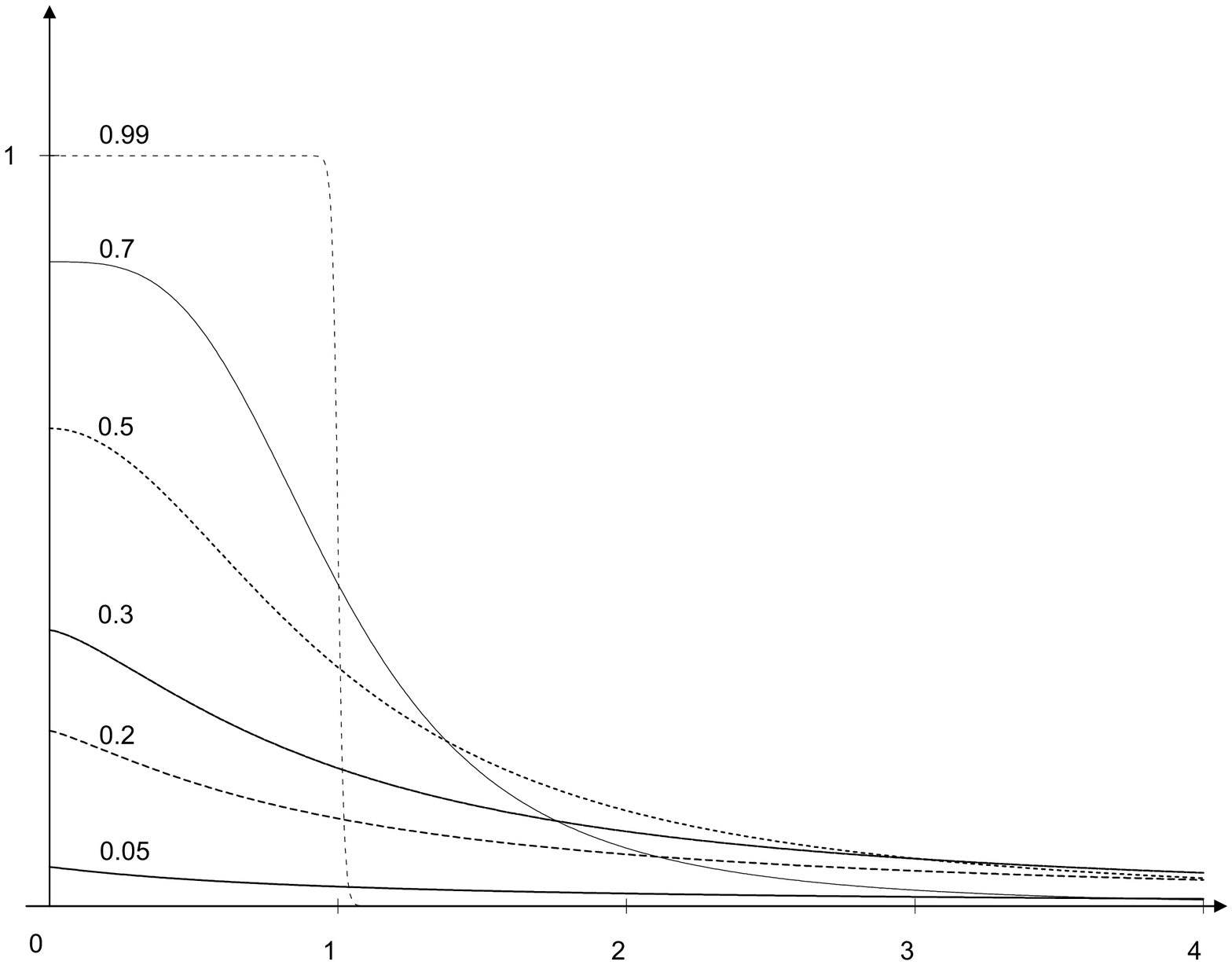}

\caption{\label{cap:Figdensity3} The densities $f_{\delta_{\alpha}}$ of
the limiting distribution of the normalized Kac process for different
values of $\alpha\in\left(0,1\right)$. The extreme case is given
by $\delta_{0}=\infty$. }
\end{figure}

\item If $0<\alpha\leq1,$ then we have\[
\Theta_{n}\;\stackrel{\mathcal{L\left(\mu\right)}}{\Longrightarrow}\;\theta_{\alpha},\]
where $\chi_{\alpha}$ denotes the random variable on $\left[0,1\right]$
with density \[
f_{\theta_{\alpha}}\left(x\right)=\frac{\sin\pi\alpha}{\pi\alpha}\frac{1}{\left(1-x^{1/\alpha}\right)^{\alpha}},\qquad\alpha\in\left(0,1\right),\]
and $\theta_{1}=1$ (cf. Fig. \ref{cap:Figdensity2}). %
\begin{figure}
\psfrag{0.98}{{\small $\alpha=0.98$}}\psfrag{0.8}{{\small $\alpha=0.8$}}

\psfrag{0.5}{{\small $\alpha=0.5$}}

\psfrag{0.3}{{\small $\alpha=0.3$}}

\psfrag{0.15}{{\small $\alpha=0.15$}}\psfrag{0.005}{{\small $\alpha=0.005$}} \psfrag{0}{{\small $0$}}\psfrag{1}{{\small $1$}}\psfrag{1.5}{{\small $1$}}\psfrag{1.5}{{\small $1.5$}}\psfrag{2}{{\small $2$}}\psfrag{2.5}{{\small $2.5$}}\psfrag{3}{{\small $3$}}\psfrag{4}{{\small $4$}}

\includegraphics[%
  width=0.60\textwidth]{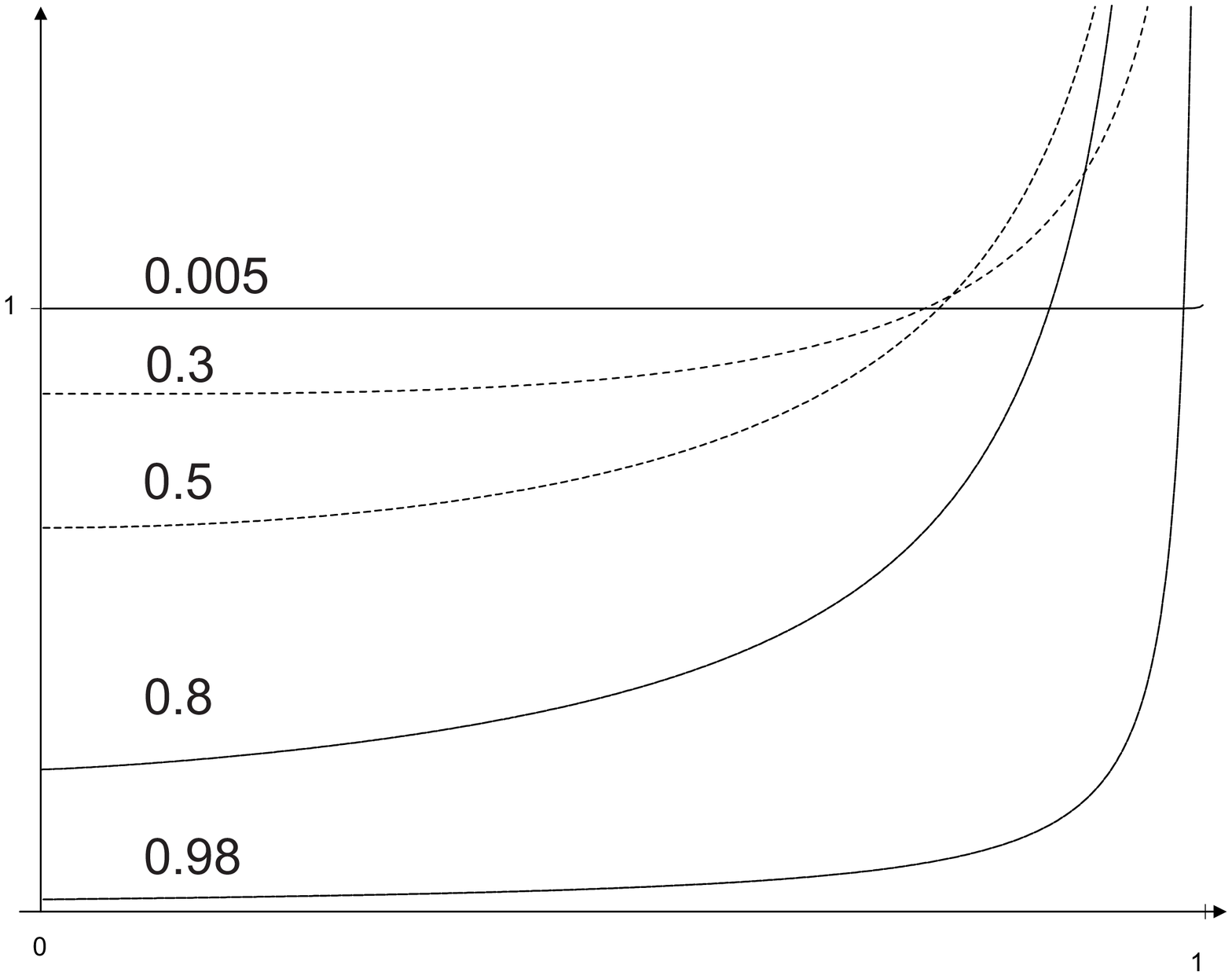}

\caption{\label{cap:Figdensity2} The densities $f_{\theta_{\alpha}}$ of
the limiting distribution of the normalized Kac process for different
values of $\alpha\in\left(0,1\right)$. The extreme case is given
by $\theta_{1}=1$. }
\end{figure}

\end{enumerate}
\end{cor}
\begin{rem*}
For $\alpha\in\left(0,1\right)$ we have \[
\lambda_{\alpha}\stackrel{\textrm{dist.}}{=}\left(\eta_{\alpha}\right)^{1-\alpha},\,\gamma_{\alpha}\stackrel{\textrm{dist.}}{=}\left(\eta_{\alpha}\right)^{-\alpha},\;\delta_{\alpha}\stackrel{\textrm{dist.}}{=}\left(\varphi_{\alpha}-1\right)^{1-\alpha},\;\textrm{and}\;\theta_{\alpha}\stackrel{\textrm{dist.}}{=}\left(\varphi_{\alpha}\right)^{-\alpha}.\]
Note, that in particular $\theta_{\frac{1}{2}}$ obeys the \emph{arc-sine
law}, i.e. it has density\[
f_{\theta_{\frac{1}{2}}}\left(x\right)=\frac{2}{\pi}\frac{1}{\sqrt{1-x^{2}}},\qquad0<x<1,\]
and $\delta_{\frac{1}{2}}$ obeys the \emph{Cauchy law}, i.e. it has
density\[
f_{\delta_{\frac{1}{2}}}\left(x\right)=\frac{2}{\pi}\frac{1}{1+x^{2}},\qquad x>0.\]

\end{rem*}
The following two theorems treat the four cases with $\alpha\in\left\{ 0,1\right\} $
not covered by Corollary \ref{coro 1}.

\begin{thm}
\label{thm:neu1} Let $A\in\mathcal{A}$ with $0<\mu\left(A\right)<\infty$
be a uniform set. If the wandering rate $\left(W_{n}\right)$ is regularly
varying with exponent $1$ then we have\begin{equation}
\Gamma_{n}\stackrel{\mathcal{L\left(\mu\right)}}{\;\Longrightarrow}\;\boldsymbol{U}\,\textrm{ and }\,\Theta_{n}\stackrel{\mathcal{L\left(\mu\right)}}{\;\Longrightarrow}\;\boldsymbol{U},\label{(result 1)}\end{equation}
where $\boldsymbol{U}$ denotes the random variable uniformly distributed
on the unit interval. 
\end{thm}
\begin{example*}
Let $f\left(0\right)=0,\; f\left(x\right)=x+x^{2}e^{-\frac{1}{x}},\quad x>0,$
and let $a\in\left(0,1\right)$ be determined by $f\left(a\right)=1.$
Define $T:\left[0,1\right]\longrightarrow\left[0,1\right]$ by\[
T\left(x\right):=\left\{ \begin{array}{ll}
f\left(x\right), & x\in\left[0,a\right],\\
\frac{x-a}{1-a}, & x\in\left(a,1\right].\end{array}\right.\]
Then the map $T$ satisfies Thaler's conditions (T1)--(T4) in \cite{thaler:95}.
Any set $A\in\mathcal{B}_{\left[0,1\right]}$ with $\lambda\left(A\right)>0$
which is bounded away from the indifferent fixed points is a uniform
set for $T.$ Furthermore, we have \[
W_{n}\sim const\cdot\frac{n}{\log\left(n\right)}\qquad\left(n\to\infty\right).\]
 Hence,\[
\frac{\log\left(n\right)}{\log\left(V_{n}\right)}\stackrel{\mathcal{L\left(\mu\right)}}{\;\Longrightarrow}\;\boldsymbol{U}\quad\textrm{and }\:\frac{\log\left(n\right)}{\log\left(Y_{n}\right)}\stackrel{\mathcal{L\left(\mu\right)}}{\;\Longrightarrow}\;\boldsymbol{U}.\]

\end{example*}
\begin{thm}
\label{theo2} Let $A\in\mathcal{A}$ with $0<\mu\left(A\right)<\infty$
be a uniformly returning set$.$ If the wandering rate $\left(W_{n}\right)$
is slowly varying, then we have\[
\Lambda_{n}\stackrel{\mathcal{L\left(\mu\right)}}{\;\Longrightarrow}\;\boldsymbol{U}\;\textrm{and }\:\Delta_{n}\stackrel{\mathcal{L\left(\mu\right)}}{\;\Longrightarrow}\;\boldsymbol{U},\]
where $\boldsymbol{U}$ denotes the random variable uniformly distributed
on the unit interval. 
\end{thm}
\begin{example*}
We consider the \emph{Lasota--Yorke} map $T:\left[0,1\right]\longrightarrow\left[0,1\right]$,
defined by\[
T\left(x\right):=\left\{ \begin{array}{ll}
\frac{x}{1-x}, & x\in\left[0,\frac{1}{2}\right],\\
2x-1, & x\in\left(\frac{1}{2},1\right].\end{array}\right.\]
This map satisfies the Thaler's conditions (i)--(iv) in \cite{Thaler:00}.
Any compact subset $A$ of $\left(0,1\right]$ with $\lambda\left(A\right)>0$
is a uniformly returning set and we have\[
W_{n}\sim\log\left(n\right)\quad\textrm{as}\quad n\to\infty.\]
Hence,\[
\frac{\log\left(Y_{n}-n\right)}{\log\left(n\right)}\stackrel{\mathcal{L\left(\mu\right)}}{\;\Longrightarrow}\;\boldsymbol{U}\quad\textrm{and}\quad\frac{\log\left(V_{n}\right)}{\log\left(n\right)}\stackrel{\mathcal{L\left(\mu\right)}}{\;\Longrightarrow}\;\boldsymbol{U}.\]

\end{example*}
Another application of the above theorem will be given in the last
section on continued fractions.

\begin{rem*}
The processes considered in \cite{KesseboehmerSlassi:05}, Theorem
1.5 and Theorem1.6, can be expressed in terms of $F$ and $G$ by
$G\left(n\right)/G\left(Z_{n}\right)$ for $\alpha=0$ and $F\left(n-Z_{n}\right)/F\left(n\right)$
for $\alpha=1$, respectively. Hence, taking the earlier result from
\cite{KesseboehmerSlassi:05} into account we have developed non-degenerated
results for all critical cases for the processes $Z_{n}$, $Y_{n}$,
and $V_{n}$.
\end{rem*}
Finally, we give a commen large deviation asymptotic for the two processes
$Z_{n}$ and $Y_{n}$, as well as a large deviation asymptotic for
the process $V_{n}$. An application of this theorem will be given
in the last section on continued fractions.

\begin{thm}
\label{thm:neue2}Let $A$ be both uniformly returning and uniform
for $f\in\mathcal{P_{\mu}}$ and let $\nu$ denote the probability
measure with density $f$. We suppose that $\left(\mu\left(A\cap\left\{ \varphi>n\right\} \right)\right)$
satisfies the following asymptotic\[
\mu\left(A\cap\left\{ \varphi>n\right\} \right)\sim n^{-1}L\left(n\right)\qquad\textrm{as}\: n\to\infty,\]
where $L$ a slowly varying function. 
\begin{itemize}
\item For $0\leq x<1$ and $y\geq0$ with $x+y\neq0$ we have\[
\nu\left(\frac{n-Z_{n}}{n}\geq x,\frac{Y_{n}-n}{n}>y\right)\sim\log\left(\frac{1+y}{x+y}\right)\cdot\frac{L\left(n\right)}{W_{n}}\qquad\textrm{as}\: n\to\infty.\]

\item For $x>0$ we have\[
\nu\left(\frac{V_{n}}{n}>x\right)\sim H\left(x\right)\frac{L\left(n\right)}{W_{n}}\qquad\textrm{as}\: n\to\infty,\]
where \[
H\left(x\right):=\left\{ \begin{array}{ll}
1-\log\left(x\right) & \quad\textrm{for }\: x\in\left(0,1\right),\\
1/x & \quad\textrm{for }\: x\geq1.\end{array}\right.\]
 
\end{itemize}
\end{thm}

\section{Proofs\label{sec:Proof-of-main}}

\subsection{Some facts from regular variation\label{sub:Classical-results-on}}

We recall the concepts of regularly varying functions and sequences
(see also \cite{BinghamGoldieTeugels:89} for a comprehensive account).
Throughout we use the convention that for two sequences $\left(a_{n}\right)$,
$\left(b_{n}\right)$ we write $a_{n}=o\left(b_{n}\right)$ if $b_{n}\not=0$
fails only for finitely many $n$ and $\lim_{n\to\infty}a_{n}/b_{n}=0$.

A measurable function $R:\mathbb{R}^{+}\rightarrow\mathbb{R}$ with
$R>0$ on $\left(a,\infty\right)$ for some $a>0$ is called \emph{regularly
varying} at $\infty$ with exponent $\rho\in\mathbb{R}$ if\[
\lim_{t\to\infty}\frac{R\left(\lambda t\right)}{R\left(t\right)}=\lambda^{\rho}\quad{\rm for\; all}\;\lambda>0.\]

A regularly varying function $L$ with exponent $\rho=0$ is called
\emph{slowly varying} at $\infty$, i.e.\[
\lim_{t\to\infty}\frac{L\left(\lambda t\right)}{L\left(t\right)}=1\quad{\rm for\; all}\;\lambda>0.\]
 Clearly, a function $R:\mathbb{R}^{+}\rightarrow\mathbb{R}$ is  regularly
varying \emph{at} $\infty$ with exponent $\rho$$\in\mathbb{R}$
if and only if\[
R\left(t\right)=t^{\rho}L\left(t\right),\quad t\in\mathbb{R}^{+},\]

for $L$ slowly varying at $\infty.$

A function $R$ is said to be \emph{regularly varying at $0$} if
$t\mapsto R\left(t^{-1}\right)$ is regularly varying at $\infty.$

A \emph{sequence} $\left(u_{n}\right)$ is \emph{regularly varying
with exponent} $\rho$ if $u_{n}=R\left(n\right)$,~$n\geq1$, for
$R:\mathbb{R}^{+}\rightarrow\mathbb{R}$ regularly varying at $\infty$
with exponent $\rho$.

The following facts will be needed in the proofs of the preparatory
lemmas and propositions of this sections, as well as for the main
theorems. 

\begin{itemize}
\item [(KL)] \textbf{Karamata's Lemma} (\cite{Feller:71}, \cite{Karamata33b}).
If $\left(a_{n}\right)$ is a regularly varying sequence with exponent
$\rho$ and if $p\geq-\rho-1,$ then \[
\lim_{n\to\infty}\frac{n^{p+1}a_{n}}{\sum_{k\leq n}k^{p}a_{k}}=p+\rho+1.\]

\item [(UA)] \textbf{Uniform asymptotics} (\cite{Seneta:76}) Let $\left(p_{n}\right)$
and $\left(q_{n}\right)$ be two positive sequences with $p_{n}\rightarrow\infty$
and $\frac{p_{n}}{q_{n}}\in\left[1/K,K\right],\; K\geq1$ for $n$
large enough. Then for every slowly varying function $L$ we have\[
\lim_{n\to\infty}\frac{L\left(p_{n}\right)}{L\left(q_{n}\right)}=1.\]

\item [(EL)] \textbf{Erickson Lemma} (\cite{Erickson:70}) Let $L\nearrow\infty$
be a monotone increasing continuous slowly varying function. Let $a_{t}\left(x\right)$
be defined by $a_{t}\left(x\right):=L^{-1}\left(xL\left(t\right)\right)$
with $x\in\left(0,\infty\right),$ where $L^{-1}\left(\cdot\right)$
denoting the inverse function of $L\left(\cdot\right)$. Then we have
for every fixed $x\in\left(0,\infty\right)$\[
a_{t}\left(x\right)\longrightarrow\infty\qquad\left(t\to\infty\right)\]
and for $0<x<y$\[
a_{t}\left(x\right)/a_{t}\left(y\right)\longrightarrow0\qquad\left(t\to\infty\right).\]

\end{itemize}

\subsection{Compactness results}

Under the assumption that $T$ is a nonsingular ergodic transformation
on $\left(X,\mathcal{A},\mu\right)$ the compactness theorem in \cite{Aaronson:97},
Section 3.6, gives the following implication. 

\begin{itemize}
\item If $R_{n}\circ T-R_{n}\stackrel{\mu}{\longrightarrow}0$ and $R_{n}\stackrel{\nu}{\Longrightarrow}R$
for some $\nu\in\mathcal{P_{\mu}}$ then $R_{n}\stackrel{\mathcal{L\left(\mu\right)}}{\Longrightarrow}R$. 
\end{itemize}
Hence, before proving the main theorems we show the following two
lemmas.

\begin{lem}
\label{lem3} Let $A\in\mathcal{A}$ be a set of positive finite measure
$\mu\left(A\right)$ and $\widetilde{L}\left(t\right)\to\infty$,
$t\to\infty$, be a slowly varying function such that\[
\widetilde{L}\left(x\right)=C\exp\left(\int_{B}^{x}\frac{\zeta\left(t\right)}{t}\; dt\right)\quad\textrm{for\; all}\; x\geq B,\]
where $C\in\left(0,\infty\right)$ and $\zeta$ a continuous function
on $\left[B,\infty\right)$ with\[
\zeta\left(x\right)\longrightarrow0\qquad\left(x\to\infty\right).\]
Then we have\begin{equation}
\frac{1}{\widetilde{L}\left(n\right)}\left(\widetilde{L}\left(Y_{n}\circ T\right)-\widetilde{L}\left(Y_{n}\right)\right)\;\stackrel{\mu}{\longrightarrow}\;0\label{eq:101}\end{equation}
and \begin{equation}
\frac{1}{\widetilde{L}\left(n\right)}\left(\widetilde{L}\left(V_{n}\circ T\right)-\widetilde{L}\left(V_{n}\right)\right)\;\stackrel{\mu}{\longrightarrow}\;0.\label{eq:102}\end{equation}

\end{lem}
\begin{proof} Without loss of generality we assume that there exists
$\delta\in\left(0,1\right)$ such that \[
\left|\zeta\left(t\right)\right|<\delta\quad\textrm{for\; all}\; t\geq B.\]
For $\varepsilon>0$ we define\[
K_{\varepsilon,n}:=\left\{ Y_{n}<\infty\;\wedge\;\frac{1}{\widetilde{L}\left(n\right)}\left|\widetilde{L}\left(Y_{n}\circ T\right)-\widetilde{L}\left(Y_{n}\right)\right|\geq\varepsilon\right\} \quad\left(n\in\mathbb{N}\right).\]
Since\begin{equation}
Y_{n}\left(T\left(x\right)\right)=\left\{ \begin{array}{ll}
Y_{n}\left(x\right)-1, & x\in\left\{ Y_{n}<\infty\right\} \cap T^{-\left(n+1\right)}A^{c},\\
n+\varphi\left(T^{n+1}\left(x\right)\right), & x\in T^{-\left(n+1\right)}A,\end{array}\right.\label{z(n)-z(n-1)}\end{equation}
we conclude\begin{eqnarray*}
K_{\varepsilon,n} & \subset & \left(\left\{ Y_{n}<\infty\right\} \cap T^{-\left(n+1\right)}A^{c}\cap\left\{ \frac{1}{\widetilde{L}\left(n\right)}\left(\widetilde{L}\left(Y_{n}\right)-\widetilde{L}\left(Y_{n}-1\right)\right)\geq\varepsilon\right\} \right)\\
 &  & \cup\left(T^{-\left(n+1\right)}A\cap\left\{ \frac{1}{\widetilde{L}\left(n\right)}\left(\widetilde{L}\left(n+\varphi\left(T^{n+1}\left(\omega\right)\right)\right)-\widetilde{L}\left(n+1\right)\right)\geq\varepsilon\right\} \right).\end{eqnarray*}
 For $n\geq B$ large enough such that $C\delta B^{-\delta}\left(n-1\right)^{\delta-1}L\left(n\right)^{-1}<\varepsilon$
we have \\
${\displaystyle \widetilde{L}\left(Y_{n}\left(\omega\right)\right)-\widetilde{L}\left(Y_{n}\left(\omega\right)-1\right)=C\exp\left(\int_{B}^{Y_{n}\left(\omega\right)-1}\frac{\zeta\left(t\right)}{t}\; dt\right)\times}$\[
\;\quad\quad\quad\quad\quad\;\quad\left[\exp\left(\int_{Y_{n}\left(\omega\right)-1}^{Y_{n}\left(\omega\right)}\frac{\zeta\left(t\right)}{t}\; dt\right)-1\right].\]
Since $\left|\zeta\left(t\right)\right|<\delta$ on $\left[B,\infty\right)$
we have by the Mean-Value Theorem\begin{eqnarray*}
\widetilde{L}\left(Y_{n}\left(\omega\right)\right)-\widetilde{L}\left(Y_{n}\left(\omega\right)-1\right) & \leq & \frac{C}{B^{\delta}}\left(\left(Y_{n}\left(\omega\right)\right)^{\delta}-\left(Y_{n}\left(\omega\right)-1\right)^{\delta}\right)\\
 & \leq & \frac{C\delta}{B^{\delta}}\left(n-1\right)^{\delta-1}.\end{eqnarray*}
Now choose $n\geq B$ large enough such that $\frac{C\delta\left(n-1\right)^{\delta-1}}{B^{\delta}\widetilde{L}\left(n\right)}<\varepsilon$.
This implies \textbf{}\[
K_{\varepsilon,n}\subset\left(T^{-\left(n+1\right)}A\cap\left\{ \frac{1}{\widetilde{L}\left(n\right)}\left(\widetilde{L}\left(n+\varphi\left(T^{n+1}\left(\omega\right)\right)\right)-\widetilde{L}\left(n+1\right)\right)\geq\varepsilon\right\} \right).\]
Similarly as above, we obtain for sufficiently large $n$ \\
${\displaystyle \widetilde{L}\left(n+\varphi\left(T^{n+1}\left(\omega\right)\right)\right)-\widetilde{L}\left(n+1\right)=C\exp\left(\int_{B}^{n+1}\frac{\zeta\left(t\right)}{t}\; dt\right)\times}$\[
\;\quad\quad\quad\quad\quad\;\quad\left[\exp\left(\int_{n+1}^{n+\varphi\left(T^{n+1}\left(\omega\right)\right)}\frac{\zeta\left(t\right)}{t}\; dt\right)-1\right].\]
 Since $\left|\zeta\left(t\right)\right|<\delta$ on $\left[B,\infty\right)$,
there exists a constant $C_{\delta}$, such that\begin{eqnarray*}
\widetilde{L}\left(n+\varphi\left(T^{n+1}\left(\omega\right)\right)\right)-\widetilde{L}\left(n+1\right) & \leq & C_{\delta}\left(\left(n+\varphi\left(T^{n+1}\left(\omega\right)\right)\right)^{\delta}-\left(n+1\right)^{\delta}\right)=:E.\end{eqnarray*}
By the Mean-Value Theorem, we have 

\[
E\leq\delta C_{\delta}n^{\delta-1}\left(\varphi\left(T^{n+1}\left(\omega\right)\right)-1\right).\]
Hence,\\
${\displaystyle \left(T^{-\left(n+1\right)}A\cap\left\{ \frac{\widetilde{L}\left(n+\varphi\left(T^{n+1}\left(\omega\right)\right)\right)-\widetilde{L}\left(n+1\right)}{\widetilde{L}\left(n\right)}\geq\varepsilon\right\} \right)}$\[
\qquad\qquad\qquad\subset T^{-\left(n+1\right)}\left(A\cap\left\{ \varphi\geq\frac{n^{1-\delta}\widetilde{L}\left(n\right)}{\delta C_{\delta}}\varepsilon+1\right\} \right).\]
Using the invariance of $\mu$ and the fact that by choice of $\delta$
we have $n^{1-\delta}\widetilde{L}\left(n\right)\to\infty$, we obtain
\[
\mu\left(K_{\varepsilon,n}\right)\leq\mu\left(A\cap\left\{ \varphi\geq\frac{n^{1-\delta}\widetilde{L}\left(n\right)}{\delta C_{\delta}}\varepsilon+1\right\} \right)\to0\quad\textrm{for}\quad n\to\infty.\]
This implies \[
\lim_{n\to\infty}\nu\left(K_{\varepsilon,n}\right)=0\quad\textrm{for all}\;\nu\in\mathcal{P}_{\mu}.\]
Using this and the fact $\left\{ Y_{n}<\infty\right\} =X$ modulo
a set of $\mu$--measure $0$, we finally conclude for all $\nu\in\mathcal{P}_{\mu}$
\[
\lim_{n\to\infty}\nu\left(\left\{ \frac{\left|\widetilde{L}\left(Y_{n}\circ T\right)-\widetilde{L}\left(Y_{n}\right)\right|}{\widetilde{L}\left(n\right)}\geq\varepsilon\right\} \right)=0.\]

The second assertion follows analogously by using the first part of
the lemma and \cite{KesseboehmerSlassi:05}, Lemma 3.1. \end{proof}

\begin{lem}
\label{lem:kompakt}Let $A\in\mathcal{A}$ be a set of positive finite
measure $\mu\left(A\right),$ then \[
\Delta_{n}\circ T-\Delta_{n}\;\stackrel{\mu}{\longrightarrow}\;0\quad\textrm{and}\quad\Lambda_{n}\circ T-\Lambda_{n}\;\stackrel{\mu}{\longrightarrow}\;0.\]

\end{lem}
\begin{proof} Let $\varepsilon>0$ be given, and let\[
K_{\varepsilon,n}:=\left\{ Y_{n}<\infty\wedge\left|\Delta_{n}\circ T-\Delta_{n}\right|\geq\varepsilon\right\} .\]
choose $n$ large enough such that $\frac{\mu\left(A\right)}{W_{n}}<\varepsilon.$
By (\ref{z(n)-z(n-1)}) we have \begin{eqnarray*}
K_{\varepsilon,n} & \subset & T^{-\left(n+1\right)}A\cap\left\{ \varphi\left(T^{n+1}\left(\omega\right)\right)-1\geq\varepsilon\frac{W_{n}}{\mu\left(A\right)}\right\} \\
 & \subset & T^{-\left(n+1\right)}\left(A\cap\left\{ \varphi\geq\varepsilon\frac{W_{n}}{\mu\left(A\right)}+1\right\} \right).\end{eqnarray*}
this implies\[
\mu\left(K_{\varepsilon,n}\right)\to0\quad\textrm{as}\: n\to\infty.\]
Thus, \[
\lim_{n\to\infty}\nu\left(K_{\varepsilon,n}\right)=0\quad\textrm{for all}\;\nu\in\mathcal{P}_{\mu}.\]
From the fact that $\left\{ Y_{n}<\infty\right\} =X$ modulo a set
of $\mu$--measure $0$ the first assertion of the lemma follows. 

The second assertion follows analogously by using the first part of
the lemma and \cite{KesseboehmerSlassi:05}, Lemma 3.2. \end{proof}

\subsection{Proofs of main theorems}

\begin{proof} (First part of Theorem \ref{thm:neu1}) For $\alpha=1$
we have $W_{n}/n\sim1/L\left(n\right)$ for a slowly varying function
$L$ with $L\left(n\right)\to\infty$. Due to the Representation Theorem
for slowly varying functions (cf. \cite{Seneta:76}) there exists
a slowly varying function $\widetilde{L}$ with the same properties
as in Lemma \ref{lem3} such that $L\left(x\right)\sim\widetilde{L}\left(x\right)$
as $x\to\infty.$ Therefore, to prove the first part it suffices to
show \begin{equation}
\frac{\widetilde{L}\left(n\right)}{\widetilde{L}\left(V_{n}\right)}\stackrel{\mathcal{L\left(\mu\right)}}{\;\longrightarrow}\;\boldsymbol{U}.\label{result2}\end{equation}
Let $A$ be a uniform set for some $f\in\mathcal{P_{\mu}}.$ First,
for every $x\in\left(0,1\right]$, we have \begin{eqnarray*}
\nu\left(\frac{\widetilde{L}\left(n\right)}{\widetilde{L}\left(V_{n}\right)}<x\right) & = & \sum_{k=0}^{n}\nu\left(Y_{n}>k+\left\lfloor a_{n}\left(x^{-1}\right)\right\rfloor ,Z_{n}=k\right)\\
 & = & \int_{A}\sum_{k=0}^{n}\hat{T}^{k}\left(f\right)\cdot1_{A\cap\left\{ \varphi>\left\lfloor a_{n}\left(x^{-1}\right)\right\rfloor \right\} }\; d\mu,\end{eqnarray*}
where $\nu$ denotes the probability measure with density $f\in P_{\mu}$
and $a_{n}\left(x^{-1}\right)=\widetilde{L}^{-1}\left(x^{-1}\widetilde{L}\left(n\right)\right).$
Note, by (EL) we have\[
a_{n}\left(x^{-1}\right)\rightarrow\infty\quad\textrm{and}\quad\frac{n}{a_{n}\left(x^{-1}\right)}\rightarrow0\qquad\textrm{for}\; n\to\infty.\]
By the asymptotic in (\ref{eq:AsymptAaronson}) and (KL) we obtain
on the one hand that for $\varepsilon\in\left(0,1\right)$ and sufficiently
large $n$\begin{eqnarray*}
\nu\left(\frac{\widetilde{L}\left(n\right)}{\widetilde{L}\left(V_{n}\right)}<x\right) & \leq & \left(1+\varepsilon\right)\mu\left(A\cap\left\{ \varphi>\left\lfloor a_{n}\left(x^{-1}\right)\right\rfloor \right\} \right)\cdot\widetilde{L}\left(n\right)\\
 & \sim & \left(1+\varepsilon\right)x.\end{eqnarray*}
This implies\[
\limsup\nu\left(\frac{\widetilde{L}\left(n\right)}{\widetilde{L}\left(V_{n}\right)}<x\right)\leq x.\]
On the other hand, we similarly obtain\[
\liminf\nu\left(\frac{\widetilde{L}\left(n\right)}{\widetilde{L}\left(V_{n}\right)}<x\right)\geq x.\]
Both inequalities give \[
\nu\left(\frac{\widetilde{L}\left(n\right)}{\widetilde{L}\left(V_{n}\right)}<x\right)\to x.\]
Now let $x>1$. Then we have\[
\nu\left(\frac{\widetilde{L}\left(n\right)}{\widetilde{L}\left(V_{n}\right)}<1\right)\leq\nu\left(\frac{\widetilde{L}\left(n\right)}{\widetilde{L}\left(Y_{n}\right)}<x\right)\leq1.\]
From this it follows that\[
\nu\left(\frac{\widetilde{L}\left(n\right)}{\widetilde{L}\left(V_{n}\right)}<x\right)\to1.\]
Hence, we obtain for all $x\in\left(0,\infty\right)$\[
\nu\left(\frac{\widetilde{L}\left(n\right)}{\widetilde{L}\left(V_{n}\right)}<x\right)\to1-\max\left(1-x,0\right).\]
Using this and Lemma \ref{lem3} the convergence in (\ref{result2})
follows by the compactness theorem. Finally, since $V_{n}\rightarrow\infty$
in probability, it is clear that the slowly varying function $\widetilde{L}$
may be replaced by any function $L$with $L\left(n\right)\sim C\cdot\widetilde{L}\left(n\right),$
$C>0$, as $n\to\infty.$ This finishes the the proof of the first
part.\end{proof}

\begin{proof} (Second part of Theorem \ref{thm:neu1}) Let $f$ and
$\widetilde{L}$ be given as in the first part of the proof of this
theorem. Since\[
\left\{ Z_{n}\leq k\right\} =\left\{ Y_{k}>n\right\} \quad\textrm{for}\:1\leq k\leq n,\]
we have for every $x\in\left(0,1\right)$ \begin{eqnarray*}
\nu\left(\frac{\widetilde{L}\left(n\right)}{\widetilde{L}\left(Y_{n}\right)}<x\right) & = & \nu\left(Y_{n}>\left\lfloor a_{n}\left(x^{-1}\right)\right\rfloor \right)\\
 & = & \nu\left(Z_{\left\lfloor a_{n}\left(x^{-1}\right)\right\rfloor }\leq n\right)\\
 & = & \int_{A}\sum_{k=0}^{n}\hat{T}^{k}\left(f\right)\cdot1_{A\cap\left\{ \varphi>\left\lfloor a_{n}\left(x^{-1}\right)\right\rfloor -k\right\} }\; d\mu,\end{eqnarray*}
where $\nu$ denotes the probability measure with density $f\in P_{\mu}$
and $a_{n}\left(x^{-1}\right)=\widetilde{L}^{-1}\left(x^{-1}\widetilde{L}\left(n\right)\right).$
Using the monotonicity of the sequence $\left(1_{A\cap\left\{ \varphi>n\right\} }\right)$
we obtain by the asymptotic in (\ref{eq:AsymptAaronson}) on the one
hand that\begin{eqnarray*}
\nu\left(\frac{\widetilde{L}\left(n\right)}{\widetilde{L}\left(Y_{n}\right)}<x\right) & \leq & \int_{A}1_{A\cap\left\{ \varphi>\left\lfloor a_{n}\left(x^{-1}\right)\right\rfloor -n\right\} }\cdot\sum_{k=0}^{n}\hat{T}^{k}\left(f\right)\; d\mu\\
 & \sim & \mu\left(A\cap\left\{ \varphi>\left\lfloor a_{n}\left(x^{-1}\right)\right\rfloor -n\right\} \right)\cdot\widetilde{L}\left(n\right).\end{eqnarray*}
This together with (KL) and (UA) implies \[
\limsup\nu\left(\frac{\widetilde{L}\left(n\right)}{\widetilde{L}\left(Y_{n}\right)}<x\right)\leq x.\]
On the other hand, we derive in a similar way \begin{eqnarray*}
\nu\left(\frac{\widetilde{L}\left(n\right)}{\widetilde{L}\left(Y_{n}\right)}<x\right) & \geq & \int_{A}1_{A\cap\left\{ \varphi>\left\lfloor a_{n}\left(x^{-1}\right)\right\rfloor \right\} }\cdot\sum_{k=0}^{n}\hat{T}^{k}\left(f\right)\; d\mu\\
 & \sim & \mu\left(A\cap\left\{ \varphi>\left\lfloor a_{n}\left(x^{-1}\right)\right\rfloor \right\} \right)\cdot\widetilde{L}\left(n\right).\end{eqnarray*}
This gives the opposite inequality\[
\liminf\nu\left(\frac{\widetilde{L}\left(n\right)}{\widetilde{L}\left(Y_{n}\right)}<x\right)\geq x.\]
Hence, we obtain\[
\frac{\widetilde{L}\left(n\right)}{\widetilde{L}\left(Y_{n}\right)}\stackrel{\nu}{\;\longrightarrow}\;\boldsymbol{U}.\]
Using this and Lemma \ref{lem3} the second assertion of the theorem
follows by the compactness theorem and the fact that $Y_{n}\rightarrow\infty$
in probability. \end{proof}

The proof of Theorem \ref{theo2} splits into two parts.

\begin{proof} (First part of Theorem \ref{theo2}) Let $A$ be a uniformly
returning set for some $f\in\mathcal{P_{\mu}}.$ Let $W_{n}\sim L\left(n\right)$
as $n\to\infty,$ without loss of generality we may assume that $L$
is monotone increasing and continuous. We have for every fixed $x\in\left(0,1\right)$\begin{eqnarray*}
\nu\left(\frac{L\left(V_{n}\right)}{L\left(n\right)}>x\right) & = & \nu\left(V_{n}>\left\lfloor a_{n}\left(x\right)\right\rfloor \right)\\
 & = & \sum_{k=0}^{n}\nu\left(Y_{n}>k+\left\lfloor a_{n}\left(x\right)\right\rfloor ,Z_{n}=k\right)\\
 & = & \sum_{k=0}^{n-\left\lfloor a_{n}\left(x\right)\right\rfloor -1}\cdots+\sum_{k=n-\left\lfloor a_{n}\left(x\right)\right\rfloor }^{n}\cdots\\
 & = & \sum_{k=0}^{n-\left\lfloor a_{n}\left(x\right)\right\rfloor -1}\int_{A}1_{A\cap\left\{ \varphi>n-k\right\} }\cdot\hat{T}^{k}\left(f\right)\; d\mu\\
 &  & \qquad\qquad+\sum_{k=n-\left\lfloor a_{n}\left(x\right)\right\rfloor }^{n}\int_{A}1_{A\cap\left\{ \varphi>\left\lfloor a_{n}\left(x\right)\right\rfloor \right\} }\cdot\hat{T}^{k}\left(f\right)\; d\mu\\
 & =: & I\left(n\right)+J\left(n\right).\end{eqnarray*}
First, we have to prove that \begin{equation}
I\left(n\right)\to1-x\quad\textrm{as}\; n\to\infty.\label{eq:11}\end{equation}
In fact, we have\begin{eqnarray*}
I\left(n\right) & = & \sum_{k=0}^{n}\int_{A}1_{A\cap\left\{ \varphi>n-k\right\} }\cdot\hat{T}^{k}\left(f\right)d\mu-\negmedspace\negmedspace\!\negmedspace\sum_{k=n-\left\lfloor a_{n}\left(x\right)\right\rfloor -2}^{n}\negmedspace\int_{A}1_{A\cap\left\{ \varphi>n-k\right\} }\cdot\hat{T}^{k}\left(f\right)d\mu\\
 & =: & \nu\left(A_{n}\right)-\widetilde{I}\left(n\right).\end{eqnarray*}
 we first note, that \[
\lim_{n\to\infty}\nu\left(A_{n}\right)=1.\]
By a similar argument as in the proof of Theorem \ref{theo2} we obtain
for all $\varepsilon\in\left(0,1\right)$ and sufficiently large $n$
on the one hand that \[
\widetilde{I}\left(n\right)\leq\left(1+\varepsilon\right)^{2}\frac{1}{W_{n}}W_{\left\lfloor a_{n}\left(x\right)\right\rfloor +2}\sim\left(1+\varepsilon\right)^{2}x.\]
On the other hand we have \[
x\left(1-\varepsilon\right)\sim\left(1-\varepsilon\right)\frac{1}{W_{n}}W_{\left\lfloor a_{n}\left(x\right)\right\rfloor +2}\leq\widetilde{I}\left(n\right).\]
Both inequalities give\[
\widetilde{I}\left(n\right)\to x\quad\textrm{as}\; n\to\infty,\]
and consequently (\ref{eq:11}) holds.

Now we prove that \begin{equation}
J\left(n\right)\to0\quad\textrm{as}\; n\to\infty.\label{eq:12}\end{equation}
In fact, we have for all $\varepsilon\in\left(0,1\right)$ and $n$
sufficiently large \begin{eqnarray*}
J\left(n\right) & \leq & \left(1+\varepsilon\right)^{2}\frac{\left\lfloor a_{n}\left(x\right)\right\rfloor +1}{W_{n}}\cdot\mu\left(A\cap\left\{ \varphi>\left\lfloor a_{n}\left(x\right)\right\rfloor \right\} \right)\\
 & \sim & x\left(1+\varepsilon\right)^{2}\frac{\left\lfloor a_{n}\left(x\right)\right\rfloor \mu\left(A\cap\left\{ \varphi>\left\lfloor a_{n}\left(x\right)\right\rfloor \right\} \right)}{W_{\left\lfloor a_{n}\left(x\right)\right\rfloor }}\to0.\end{eqnarray*}
This gives (\ref{eq:12}). From this and (\ref{eq:11}) it follows
that \[
\nu\left(\frac{L\left(V_{n}\right)}{L\left(n\right)}>x\right)\to1-x\quad\textrm{as}\; n\to\infty\:\textrm{for all }\; x\in\left(0,1\right).\]
Now let $x\in\left[1,\infty\right)$. Then we have for $x\geq1$ that\begin{eqnarray*}
\nu\left(\frac{L\left(V_{n}\right)}{L\left(n\right)}>x\right) & \leq & \nu\left(\frac{L\left(V_{n}\right)}{L\left(n\right)}>1\right)\\
 & \leq & \nu\left(\frac{V_{n}}{n}>1\right)\to0\quad\textrm{as}\; n\to\infty.\end{eqnarray*}
Combining the above we get for all $x\in\left(0,\infty\right)$\[
\nu\left(\frac{L\left(V_{n}\right)}{L\left(n\right)}\leq x\right)\to1-\max\left\{ 1-x,0\right\} \quad\textrm{as}\; n\to\infty.\]
Finally, since $V_{n}\to\infty$ in probability, it is clear that
the slowly varying function $L$ may be replaced by any function $L_{1}$with
$L_{1}\left(n\right)\sim C\cdot L\left(n\right),$ $C>0$, as $n\to\infty.$
Hence, by Lemma \ref{lem:kompakt} and the compactness result the
theorem follows. \end{proof}

\begin{proof}(Second part of Theorem \ref{theo2}) Let $f$ and $L$
be given as in the first part of the proof of this theorem. Then for
every fixed $x\in\left(0,1\right)$ we have\begin{eqnarray*}
\nu\left(\frac{L\left(Y_{n}-n\right)}{L\left(n\right)}\leq x\right) & = & \nu\left(Y_{n}\leq\left\lfloor a_{n}\left(x\right)\right\rfloor +n\right)\\
 & = & \nu\left(Z_{n+\left\lfloor a_{n}\left(x\right)\right\rfloor }>n\right)\\
 & = & \int_{A}\sum_{k=n+1}^{n+\left\lfloor a_{n}\left(x\right)\right\rfloor }1_{A\cap\left\{ \varphi>n+\left\lfloor a_{n}\left(x\right)\right\rfloor -k\right\} }\cdot\hat{T}^{k}\left(f\right)\; d\mu,\end{eqnarray*}
where $\nu$ denotes the probability measure with density $f\in P_{\mu}$
and $a_{n}\left(x\right)=L^{-1}\left(xL\left(n\right)\right).$ Note,
by (EL) we have\[
a_{n}\left(x\right)\rightarrow\infty\quad\textrm{and}\quad\frac{a_{n}\left(x\right)}{n}\rightarrow0\qquad\textrm{for}\; n\to\infty.\]
A similar arguments as in \cite{KesseboehmerSlassi:05}, Lemma 3.3,
show that for all $\varepsilon\in\left(0,1\right)$ there exists \emph{}$n_{0}$
such that for all $n\geq n_{0}$ and $k\in\left[n,n+a_{n}\left(x\right)\right]$
we have uniformly on $A$\[
\left(1-\varepsilon\right)\frac{1}{W_{n}}\leq\hat{T}^{k}\left(f\right)\leq\left(1+\varepsilon\right)^{2}\frac{1}{W_{n}}.\]
 From this it follows on the one hand that, for $n$ sufficiently
large, \[
\nu\left(\frac{L\left(Y_{n}-n\right)}{L\left(n\right)}\leq x\right)\leq\left(1+\varepsilon\right)^{2}\frac{1}{W_{n}}W_{\left[a_{n}\left(x\right)\right]-1}\sim\left(1+\varepsilon\right)^{2}x.\]
 Similarly for $n$ sufficiently large,\[
x\left(1-\varepsilon\right)\sim\left(1-\varepsilon\right)\frac{1}{W_{n}}W_{\left[a_{n}\left(x\right)\right]-1}\leq\nu\left(\frac{L\left(Y_{n}-n\right)}{L\left(n\right)}\leq x\right).\]
Combining these inequalities we get\[
x\left(1-\varepsilon\right)\leq\liminf\nu\left(\frac{L\left(Y_{n}-n\right)}{L\left(n\right)}\leq x\right)\leq\limsup\nu\left(\frac{L\left(Y_{n}-n\right)}{L\left(n\right)}\leq x\right)\leq\left(1+\varepsilon\right)^{2}x.\]
Since $\varepsilon$ was arbitrary, we conclude \begin{equation}
\nu\left(\frac{L\left(Y_{n}-n\right)}{L\left(n\right)}\leq x\right)\longrightarrow x\quad\textrm{as}\quad n\to\infty\qquad\textrm{for all}\; x\in\left(0,1\right).\label{eq:prob0}\end{equation}
Now let $x\in\left[1,\infty\right)$. Then we have\begin{eqnarray*}
\nu\left(\frac{L\left(Y_{n}-n\right)}{L\left(n\right)}>x\right) & \leq & \nu\left(\frac{L\left(Y_{n}-n\right)}{L\left(n\right)}>1\right)\\
 & \leq & \nu\left(\frac{Y_{n}}{n}-1>1\right)\to0\quad\textrm{as}\; n\to\infty.\end{eqnarray*}
From this and (\ref{eq:prob0}) it follows that \[
\nu\left(\frac{L\left(Y_{n}-n\right)}{L\left(n\right)}\leq x\right)\longrightarrow1-\max\left(1-x,0\right)\quad\textrm{as}\quad n\to\infty\qquad\textrm{for all}\; x\in\left(0,\infty\right).\]
Finally, since $Y_{n}-n\to\infty$ in probability, it is clear that
the slowly varying function $L$ may be replaced by any function $L_{1}$with
$L_{1}\left(n\right)\sim C\cdot L\left(n\right),$ $C>0$, as $n\to\infty.$
From this, Lemma \ref{lem:kompakt} and the compactness result the
theorem follows. \end{proof}

Finally, we prove the large deviation asymptotic stated in Theorem
\ref{thm:neue2}.

\begin{proof} (First part of Theorem \ref{thm:neue2}) Let $0\leq x<1$
and $y\geq0$ be fixed with $x+y\neq0$. We have \begin{eqnarray*}
\nu\left(\frac{n-Z_{n}}{n}\geq x,\frac{Y_{n}-n}{n}>y\right) & = & \nu\left(Z_{n}\leq\left\lfloor n\left(1-x\right)\right\rfloor ,Y_{n}>\left\lfloor n\left(1+y\right)\right\rfloor \right)\\
 & = & \nu\left(Z_{\left\lfloor n\left(1+y\right)\right\rfloor }\leq\left\lfloor n\left(1-x\right)\right\rfloor \right)\\
 & = & \int_{A}\sum_{k=0}^{\left\lfloor n\left(1-x\right)\right\rfloor }1_{A\cap\left\{ \varphi>\left\lfloor n\left(1+y\right)\right\rfloor -k\right\} }\cdot\hat{T}^{k}\left(f\right)\; d\mu.\end{eqnarray*}
For $\delta\in\left(0,1-x\right)$ and $\varepsilon\in\left(0,1\right)$
fixed but arbitrary we divide the above sum into two parts as follows.\[
\nu\left(\frac{n-Z_{n}}{n}\geq x,\frac{Y_{n}-n}{n}>y\right)=\sum_{k=0}^{\left\lfloor n\delta\right\rfloor -1}\cdots+\sum_{k=\left\lfloor n\delta\right\rfloor }^{\left\lfloor n\left(1-x\right)\right\rfloor }\cdots=:I\left(n\right)+J\left(n\right).\]
By monotonicity of $\left(1_{A\cap\left\{ \varphi>n\right\} }\right)$
we first have\begin{eqnarray*}
I\left(n\right) & \leq & \int_{A}1_{A\cap\left\{ \varphi>\left\lfloor n\left(1+y\right)\right\rfloor -\left\lfloor n\delta\right\rfloor +1\right\} }\cdot\sum_{k=0}^{\left\lfloor n\delta\right\rfloor -1}\hat{T}^{k}\left(f\right)\; d\mu.\end{eqnarray*}
Using (\ref{eq:AsymptAaronson}) and the fact that $A$ is uniform
for $f$, we obtain for sufficiently large $n$\begin{eqnarray*}
I\left(n\right) & \leq & \left(1+\varepsilon\right)^{2}\frac{\left\lfloor n\delta\right\rfloor -1}{\left\lfloor n\left(1+y\right)\right\rfloor -\left\lfloor n\delta\right\rfloor +1}\cdot\frac{L\left(\left\lfloor n\left(1+y\right)\right\rfloor -\left\lfloor n\delta\right\rfloor +1\right)}{W_{\left\lfloor n\delta\right\rfloor -1}}\\
 & \sim & \left(1+\varepsilon\right)^{2}\frac{\delta}{y+1-\delta}\frac{L\left(n\right)}{W_{n}}\quad\textrm{as}\; n\to\infty.\end{eqnarray*}
Thus, \[
\limsup_{n\to\infty}\frac{W_{n}}{L\left(n\right)}\cdot I\left(n\right)\leq\left(1+\varepsilon\right)^{2}\frac{\delta}{y+1-\delta}.\]
Letting $\delta\to0$, we observe\begin{equation}
I\left(n\right)=o\left(\frac{L\left(n\right)}{W_{n}}\right),\quad\textrm{as}\: n\to\infty.\label{I(n)}\end{equation}
For the second part of the sum we have to show that\begin{equation}
J\left(n\right)\sim\frac{L\left(n\right)}{W_{n}}\cdot\log\left(\frac{1+y}{x+y}\right)\qquad\textrm{as}\: n\to\infty.\label{J(n)}\end{equation}
A similarly argument as in \cite{KesseboehmerSlassi:05}, Lemma 3.3,
shows that for all $n$ sufficiently large and $k\in\left[\left\lfloor n\delta\right\rfloor ,\left\lfloor n\left(1-x\right)\right\rfloor \right]$
we have uniformly on $A$\begin{equation}
\left(1-\varepsilon\right)\frac{1}{W_{n}}\leq\hat{T}^{k}\left(f\right)\leq\left(1+\varepsilon\right)^{2}\frac{1}{W_{n}}.\label{eq:(right W_n)}\end{equation}
Hence, using the right-hand side of (\ref{eq:(right W_n)}) and (UA),
we obtain for $n$ sufficiently large\begin{eqnarray*}
J\left(n\right) & \leq & \frac{\left(1+\varepsilon\right)^{2}}{W_{n}}\cdot\sum_{k=\left\lfloor n\left(1+y\right)\right\rfloor -\left\lfloor n\left(1-x\right)\right\rfloor }^{\left\lfloor n\left(1+y\right)\right\rfloor -\left\lfloor n\delta\right\rfloor }\mu\left(A\cap\left\{ \varphi>k\right\} \right)\\
 & \sim & \left(1+\varepsilon\right)^{2}\frac{L\left(n\right)}{W_{n}}\cdot\sum_{k=\left\lfloor n\left(1+y\right)\right\rfloor -\left\lfloor n\left(1-x\right)\right\rfloor }^{\left\lfloor n\left(1+y\right)\right\rfloor -\left\lfloor n\delta\right\rfloor }\frac{1}{k}\\
 & \sim & \left(1+\varepsilon\right)^{2}\frac{L\left(n\right)}{W_{n}}\cdot\log\left(\frac{1+y-\delta}{x+y}\right)\quad\textrm{as}\; n\to\infty.\end{eqnarray*}
This implies \[
\limsup_{n\to\infty}\frac{W_{n}}{L\left(n\right)}\cdot J\left(n\right)\leq\left(1+\varepsilon\right)^{3}\log\left(\frac{1+y-\delta}{x+y}\right).\]
Similarly, using the left-hand side of (\ref{eq:(right W_n)}), we
get \[
\liminf_{n\to\infty}\frac{W_{n}}{L\left(n\right)}\cdot J\left(n\right)\geq\left(1+\varepsilon\right)^{2}\log\left(\frac{1+y-\delta}{x+y}\right).\]
Since $\varepsilon$ and $\delta$ were arbitrary, (\ref{J(n)}) holds.
Combining (\ref{I(n)}) and (\ref{J(n)}) proves then the claim of
the theorem.\qed\end{proof}

\begin{proof}(Second part of Theorem \ref{thm:neue2}) First, let
$x\in\left(0,1\right)$. We have\begin{eqnarray*}
\nu\left(\frac{V_{n}}{n}>x\right) & = & \sum_{k=0}^{n}\nu\left(Y_{n}>k+\left\lfloor nx\right\rfloor ,Z_{n}=k\right)\\
 & = & \sum_{k=0}^{n-\left\lfloor nx\right\rfloor -1}\int_{A}1_{A\cap\left\{ \varphi>n-k\right\} }\cdot\hat{T}^{k}\left(f\right)\; d\mu\\
 &  & +\sum_{k=n-\left\lfloor nx\right\rfloor }^{n}\int_{A}1_{A\cap\left\{ \varphi>\left\lfloor nx\right\rfloor \right\} }\cdot\hat{T}^{k}\left(f\right)\; d\mu\\
 & =: & I\left(n\right)+J\left(n\right).\end{eqnarray*}
Let $\delta\in\left(0,1-x\right)$ and $\varepsilon\in\left(0,1\right)$be
fixed but arbitrary. First, we prove that\begin{equation}
J\left(n\right)\sim\frac{L\left(n\right)}{W_{n}},\quad\textrm{as}\: n\to\infty.\label{(J(1))}\end{equation}
In fact, we have for sufficiently large $n$\begin{eqnarray*}
J\left(n\right) & \leq & \left(1+\varepsilon\right)^{2}\frac{\left\lfloor nx\right\rfloor -1}{W_{n}}\cdot\mu\left(A\cap\left\{ \varphi>\left\lfloor nx\right\rfloor \right\} \right)\\
 & \sim & \left(1+\varepsilon\right)^{2}\frac{L\left(n\right)}{W_{n}}\quad\textrm{as}\; n\to\infty.\end{eqnarray*}
Similarly we get\begin{eqnarray*}
J\left(n\right) & \geq & \left(1-\varepsilon\right)\frac{\left\lfloor nx\right\rfloor -1}{W_{n}}\cdot\mu\left(A\cap\left\{ \varphi>\left\lfloor nx\right\rfloor \right\} \right)\\
 & \sim & \left(1-\varepsilon\right)\frac{L_{\mu}\left(n\right)}{W_{n}}\quad\textrm{as}\; n\to\infty.\end{eqnarray*}
Combining both inequality (\ref{(J(1))}) follows. 

Now we have to prove that \begin{equation}
I\left(n\right)\sim-\log\left(x\right)\cdot\frac{L\left(n\right)}{W_{n}}\qquad\textrm{as}\: n\to\infty.\label{II(n)}\end{equation}
divide $I\left(n\right)$ into two parts as follows\[
I\left(n\right)=\sum_{k=0}^{\left\lfloor n\delta\right\rfloor -1}\cdots+\sum_{k=\left\lfloor n\delta\right\rfloor }^{n-\left\lfloor nx\right\rfloor -1}\cdots=:I_{1}\left(n\right)+I_{2}\left(n\right).\]
Using the monotonicity of $\left(1_{A\cap\left\{ \varphi>n\right\} }\right)$,
the fact that $A$ is uniformly for $f$, and (\ref{eq:AsymptAaronson})
we obtain, for $n$ sufficiently large, \begin{eqnarray*}
I_{1}\left(n\right) & \leq & \int_{A}1_{A\cap\left\{ \varphi>n-\left\lfloor n\delta\right\rfloor +1\right\} }\cdot\sum_{k=0}^{\left\lfloor n\delta\right\rfloor -1}\hat{T}^{k}\left(f\right)\; d\mu\\
 & \leq & \left(1+\varepsilon\right)^{2}\frac{\left\lfloor n\delta\right\rfloor -1}{n-\left\lfloor n\delta\right\rfloor +1}\cdot\frac{L_{\mu}\left(n-\left\lfloor n\delta\right\rfloor +1\right)}{W_{\left\lfloor n\delta\right\rfloor -1}}\\
 & \sim & \left(1+\varepsilon\right)^{2}\frac{\delta}{1-\delta}\frac{L_{\mu}\left(n\right)}{W_{n}}\quad\textrm{as}\; n\to\infty.\end{eqnarray*}
Consequently,\begin{equation}
I_{1}\left(n\right)=o\left(\frac{L_{\mu}\left(n\right)}{W_{n}}\right),\quad\textrm{as}\: n\to\infty.\label{eq:I_1(1)}\end{equation}
Now using the fact that $A$ is uniformly returning for $f$ we have,
for $n$ sufficiently large, \begin{eqnarray*}
I_{2}\left(n\right) & \leq & \frac{\left(1+\varepsilon\right)^{2}}{W_{n}}\cdot\sum_{k=\left\lfloor nx\right\rfloor +1}^{n-\left\lfloor n\delta\right\rfloor }\mu\left(A\cap\left\{ \varphi>k\right\} \right)\\
 & \sim & \left(1+\varepsilon\right)^{2}\frac{L_{\mu}\left(n\right)}{W_{n}}\cdot\log\left(\frac{1-\delta}{x}\right)\quad\textrm{as}\; n\to\infty.\end{eqnarray*}
This implies \[
\limsup_{n\to\infty}\frac{W_{n}}{L_{\mu}\left(n\right)}\cdot I_{2}\left(n\right)\leq\left(1+\varepsilon\right)^{3}\log\left(\frac{1-\delta}{x}\right).\]
Similarly, we get \[
\liminf_{n\to\infty}\frac{W_{n}}{L_{\mu}\left(n\right)}\cdot I_{2}\left(n\right)\geq\left(1+\varepsilon\right)^{2}\log\left(\frac{1-\delta}{x}\right).\]
Since $\varepsilon$ and $\delta$ were arbitrary, we have\begin{equation}
I_{2}\left(n\right)\sim-\log\left(x\right)\cdot\frac{L\left(n\right)}{W_{n}},\quad\textrm{as}\: n\to\infty.\label{eq:I_2(1)}\end{equation}
The asymptotics (\ref{eq:I_1(1)}) and (\ref{eq:I_2(1)}) prove (\ref{II(n)}).
Combining (\ref{(J(1))}) and (\ref{II(n)}) proves the second part
of the theorem for $x\in\left(0,1\right)$. 

Now we consider the case $x\geq1$. Since\begin{eqnarray*}
\nu\left(\frac{V_{n}}{n}>x\right) & = & \sum_{k=0}^{n}\int_{A}1_{A\cap\left\{ \varphi>\left\lfloor nx\right\rfloor \right\} }\cdot\hat{T}^{k}\left(f\right)\; d\mu\end{eqnarray*}
we have, for $n$ sufficiently large, \begin{eqnarray*}
\nu\left(\frac{V_{n}}{n}>x\right) & \leq & \mu\left(A\cap\left\{ \varphi>\left\lfloor nx\right\rfloor \right\} \right)\frac{n}{W_{n}}\left(1+\varepsilon\right)\\
 & \sim & \frac{L\left(n\right)}{xW_{n}}\left(1+\varepsilon\right).\end{eqnarray*}
Similarly, we obtain the reverse inequality proving the statement
in the theorem for $x\geq1$.\qed\end{proof}

\section{Application to continued fraction}

Any irrational number $x\in\I:=\left[0,1\right]\setminus\Q$ has a
simple infinite continued fraction expansion\[
x=\frac{1}{\kappa_{1}\left(x\right)+{\displaystyle \frac{1}{\kappa_{2}\left(x\right)+\cdots}}},\]
where the unique \emph{continued fraction digits} $\kappa_{n}\left(x\right)$
are from the positive integers $\mathbb{N}$. The \emph{Gauss transformation
$G:\I\rightarrow\I$} is given by\[
G(x):=\frac{1}{x}-\left\lfloor \frac{1}{x}\right\rfloor ,\]
where $\left\lfloor x\right\rfloor $ denotes the greatest integer
not exceeding $x\in\mathbb{R}.$ Write $G^{n}$ for the $n$-th iterate
of $G,$ $n\in\mathbb{N}_{0}=\left\{ 0,1,2,\ldots\right\} $ with
$G^{0}=\textrm{id.}$ It is then well known that for all $n\in\mathbb{N}$,
we have \[
\kappa_{n}(x)=\left\lfloor \frac{1}{G^{n-1}x}\right\rfloor .\]
 Clearly, the $\kappa_{n},$ $n\in\mathbb{N}$, define random variables
on the measure space $\left(\I,\mathcal{B},\mathbb{P}\right)$, where
$\mathcal{B}$ denotes the Borel $\sigma$-algebra of $\I$ and $\mathbb{P}$
some probability measure on $\mathcal{B}$. Then each $\kappa_{n}$
has infinite expectation with respect to the Lebesgue measure on $[0,1]$,
which we will denote by $\lambda$.

Given $n\geq1,$ we define the Process\[
\psi_{n}\left(x\right):=\max\left\{ p\in\mathbb{N}_{0}\;:\;\sum_{i=1}^{p}\kappa_{i}\left(x\right)\leq n\right\} ,\quad x\in\I,\]
and we concider the Process\begin{equation}
\sigma_{n}\left(x\right):=\kappa_{\psi_{n}\left(x\right)+1},\quad x\in\I.\label{eq:process}\end{equation}

In this paper we want to demonstrate how infinite ergodic theory can
be employed to derive new insights into the stochastic structure of
the Process $\left(\sigma_{n}\right)$. The underlying dynamical system
will be given by the Farey map. 

This process turns out to be related to the total waiting time processes
considered in the first part of this paper. This allows us to derive
the following main theorem. Its proof will be postponed to the end
of Subsection \ref{sub:DistributionalProof}.

\begin{thm}
\label{thm:main} Let $\sigma_{n}$ be the process given in (\ref{eq:process}).
Then the following holds.
\begin{enumerate}
\item We have \begin{equation}
\frac{\log\left(\sigma_{n}\right)}{\log\left(n\right)}\stackrel{\mathcal{L\left(\mu\right)}}{\longrightarrow}\boldsymbol{U},\label{gl3}\end{equation}
 where the random variable $\boldsymbol{U}$ is uniformly distributed
on the unit interval. 
\item For any $\nu\in\mathcal{D}$ and $x\in\left(0,1\right)$ we have \begin{equation}
\nu\left(\frac{\sigma_{n}}{n}>x\right)\sim\frac{H\left(x\right)}{\log\left(n\right)}\quad\textrm{as}\quad n\to\infty,\label{gl2}\end{equation}
where \[
H\left(x\right):=\left\{ \begin{array}{ll}
1-\log\left(x\right) & \quad\textrm{for }\: x\in\left(0,1\right),\\
1/x & \quad\textrm{for }\: x\geq1.\end{array}\right.\]

\end{enumerate}
\end{thm}

\subsection{Farey vs. Gauss map\label{sub:Farey-vs.-Gauss}}

We consider the Farey map $T:\left[0,1\right]\rightarrow\left[0,1\right],$
defined by\[
T\left(x\right):=\left\{ \begin{array}{ll}
T_{0}\left(x\right), & x\in\left[0,\frac{1}{2}\right],\\
T_{1}\left(x\right), & x\in\left(\frac{1}{2},1\right],\end{array}\right.\]
where\[
T_{0}\left(x\right):=\frac{x}{1-x}\qquad\textrm{and}\qquad T_{1}\left(x\right):=\frac{1}{x}-1.\]

It is known that $\left(\left[0,1\right],T,\mathcal{B},\mu\right)$
defines a conservative ergodic measure preserving dynamical system,
where $\mu$ denotes the $\sigma$- finite invariant measure with
density $h\left(x\right):=\frac{d\mu}{d\lambda}\left(x\right)=\frac{1}{x}$
. Also any Borel set $A\in\mathcal{B}$ with $\lambda\left(A\right)>0$
which is bounded away from the indifferent fixed point $0$ is a uniform
set. Furthermore, from \cite{KessboehmerSlassi:05}, Lemma 3.3, we
know that the set $K_{1}:=\left(\frac{1}{2},1\right]$ is uniformly
returning for any $f\in\mathcal{D}$, where \[
\mathcal{D}:=\left\{ f\in\mathcal{P}_{\mu}:\; f\in\mathcal{C}^{2}\left(\left(0,1\right)\right)\;\textrm{with}\; f'>0\;\textrm{and}\; f''\leq0\right\} .\]
For the wandering rate we have \[
W_{n}:=W_{n}\left(K_{1}\right)=\int_{\frac{1}{n+2}}^{1}\frac{1}{x}\, dx=\log\left(n+2\right)\sim\log\left(n\right)\qquad\left(n\to\infty\right).\]

The inverse branches of the Farey map are \begin{eqnarray*}
u_{0}\left(x\right) & := & \left(T_{0}\right)^{-1}\left(x\right)={\displaystyle \frac{x}{1+x}},\\
u_{1}\left(x\right) & := & \left(T_{1}\right)^{-1}\left(x\right)={\displaystyle \frac{1}{1+x}}.\end{eqnarray*}
For $x\neq0$ the map $u_{0}\left(x\right)$ is conjugated to the
right translation $x\mapsto F\left(x\right):=x+1,$ i.e.\[
u_{0}=J\circ F\circ J\qquad\textrm{with}\quad J\left(x\right)=J^{-1}\left(x\right)=\frac{1}{x}.\]
This shows that for the $n$-th iterate we have\begin{equation}
u_{0}^{n}\left(x\right)=J\circ F^{n}\circ J\left(x\right)=\frac{x}{1+nx}.\label{u0}\end{equation}
Moreover, we have $u_{1}\left(x\right)=J\circ F\left(x\right).$ 

Let $\mathcal{F}=\left\{ K_{n}\right\} _{n\geq1}$ be the countable
collection of pairwise disjoint subintervals of $\left[0,1\right]$
given by $K_{n}:=\left(\frac{1}{n+1},\frac{1}{n}\right]$. Setting
$A_{0}=\left[0,1\right)$, it is easy to check that $T\left(K_{n}\right)=K_{n-1}$
for all $n\geq1.$ The \emph{first entry time} $e:\I\rightarrow\mathbb{N}$
into the interval $K_{1}$ is defined as\[
e\left(x\right):=\min\left\{ k\geq0:\; T^{k}\left(x\right)\in K_{1}\right\} .\]
Then the first entry time is connected to the first digit in the continued
fraction expansion by\[
\kappa_{1}\left(x\right)=1+e\left(x\right)\quad\textrm{and}\quad\varphi\left(x\right)=\kappa_{1}\circ T\left(x\right),\quad x\in\I.\]
We now consider the \emph{induced map} $S:\I\rightarrow\I$ defined
by\[
S\left(x\right):=T^{e\left(x\right)+1}\left(x\right).\]
Since for all $n\ge1$ \[
\left\{ x\in\I:\; e\left(x\right)=n-1\right\} =K_{n}\cap\I,\]
we have by (\ref{u0}) for any $x\in K_{n}\cap\I$ \[
S\left(x\right)=T^{n}\left(x\right)=T_{1}\circ T_{0}^{n-1}\left(x\right)=\frac{1}{x}-n=\frac{1}{x}-\kappa_{1}(x).\]
This implies that the induced transformation $S$ coincides with Gauss
map $G$ on $\I$.

\subsection{Renewal theory for continued fractions \label{sub:DistributionalProof}}

In the next lemma we connect the number theoretical process $\sigma_{n}$
defined in (\ref{eq:process}) with the total waiting time process
$V_{n}$ defined with respect to the Farey map.

Let $\left(\tau_{n}\right)_{n\in\mathbb{N}}$ be the sequence of return
times, i.e. integer valued positive random variables defined recursively
by \begin{eqnarray*}
\tau_{1}\left(x\right) & := & \varphi(x)=\inf\{ p\geq1:\; T^{p}(x)\in K_{1}\},\quad x\in X,\\
\tau_{n}\left(x\right) & := & \inf\{ p\geq1:\; T^{p+\sum_{k=1}^{n-1}\tau_{k}\left(x\right)}(x)\in K_{1}\},\quad x\in X.\end{eqnarray*}

The \emph{renewal process} is then given by\[
N_{n}(x):=\left\{ \begin{array}{ll}
\max\{ k\leq n\;:\; S_{k}\left(x\right)\leq n\}, & x\in A_{n}=\bigcup_{k=0}^{n}T^{-k}K_{1},\\
0, & \textrm{else,}\end{array}\right.\]
where\[
S_{0}:=0,\qquad S_{n}:=\sum_{k=1}^{n}\tau_{k},\quad n\in\mathbb{N}.\]

\begin{lem}
\label{hauptlemma} for all $x\in\I$ and $n\geq1$ we have Let $K_{1}:=\left(\frac{1}{2},1\right]$
and $A_{n}:=\bigcup_{k=0}^{n}T^{-k}K_{1}$. Then for the process $\sigma_{n}$
defined in (\ref{eq:process}) we have for all $x\in\I$ and $n\geq1$\[
\sigma_{n}\left(x\right)=\left\{ \begin{array}{ll}
V_{n-1}\left(x\right),\qquad & x\in A_{n-1},\\
1+Y_{n-1}\left(x\right), & \textrm{else.}\end{array}\right.\]

\end{lem}
\begin{proof}As a consequence of the observations in Subsection \ref{sub:Farey-vs.-Gauss}
we will argue as follows. For $x\in\I\cap A_{n-1}^{c}$ we have that
$\kappa_{1}\left(x\right)>n$ implies $\psi_{n}\left(x\right)=0$.
For $x\in\I\cap A_{n-1}$ we distinguish two cases. Either the process
starts in $x\in K_{1}$, then we have $\kappa_{1}\left(x\right)=1$
and inductively for $n\ge2$ \[
\kappa_{n}\left(x\right)=\tau_{n-1}\left(x\right),\]
 or the process starts in $x\in K_{1}^{c}$, then we have $\kappa_{1}\left(x\right)=1+\tau_{1}\left(x\right)$
and inductively for $n\ge2$\[
\kappa_{n}\left(x\right)=\tau_{n}\left(x\right).\]
This implies that\[
\psi_{n}\left(x\right)=\left\{ \begin{array}{ll}
N_{n-1}\left(x\right)+1,\qquad & x\in K_{1},\\
N_{n-1}\left(x\right), & \textrm{else.}\end{array}\right.\]
Hence, we have for $x\in\I\cap A_{n-1}$\[
\kappa_{\psi_{n}\left(x\right)+1}=\tau_{N_{n-1}\left(x\right)+1}\]
and for $x\in\I\cap A_{n-1}^{c}$\[
\kappa_{\psi_{n}\left(x\right)+1}=\kappa_{1}\left(x\right)=1+\tau_{1}\left(x\right).\]
From this the assertion follows.  \end{proof}

After these preparations we are now in the position to give the proof
of Theorem \ref{thm:main}.

\begin{proof}(Theorem \ref{thm:main}) The asymptotic (\ref{gl3})
is an immediate consequence of Lemma \ref{hauptlemma}, the second
part of Theorem \ref{theo2}, and the fact that $W_{n}\sim\log\left(n\right)$.

Finally, (\ref{gl2}) follows from the second part of Theorem \ref{thm:neue2}
by observing that $K_{1}:=\left(\frac{1}{2},1\right]$ is uniformly
returning for any $f\in\mathcal{D}$ and that \[
\mu\left(A_{1}\cap\left\{ \varphi>n\right\} \right)=\int_{\frac{n+2}{n+3}}^{1}\frac{1}{x}\: dx\;\sim\;\frac{1}{n}\quad\textrm{as}\quad n\to\infty.\]
 \end{proof}


\begin{thebibliography}{BGT89}

\bibitem[Aar97]{Aaronson:97}
J.~Aaronson.
\newblock {\em An introduction to infinite ergodic theory}, volume~50 of {\em
  Mathematical Surveys and Monographs}.
\newblock American Mathematical Society, Providence, RI, 1997.

\bibitem[BGT89]{BinghamGoldieTeugels:89}
N.~H. Bingham, C.~M. Goldie, and J.~L. Teugels.
\newblock {\em Regular variation}, volume~27 of {\em Encyclopedia of
  Mathematics and its Applications}.
\newblock Cambridge University Press, Cambridge, 1989.

\bibitem[Eri70]{Erickson:70}
K.~B. Erickson.
\newblock Strong renewal theorems with infinite mean.
\newblock {\em Trans. Amer. Math. Soc.}, 151:263--291, 1970.

\bibitem[Fel71]{Feller:71}
W.~Feller.
\newblock {\em An introduction to probability theory and its applications.
  {V}ol. {II}.}
\newblock Second edition. John Wiley \& Sons Inc., New York, 1971.

\bibitem[Kar33]{Karamata33b}
J.~Karamata.
\newblock {Sur un mode de croissance r{\'e}guli{\`e}re. Th{\'e}or{\`e}mes
  fondamentaux.}
\newblock {\em Bull. Soc. Math. Fr.}, 61:55--62, 1933.

\bibitem[KS05a]{KessboehmerSlassi:05}
M.~Kesseb{\"o}hmer and M.~Slassi.
\newblock A distributional limit law for continued fraction digit sums.
\newblock {\em arXiv:math. NT/0509559}, pages 1--15, 2005.

\bibitem[KS05b]{KesseboehmerSlassi:05}
M.~Kesseb{\"o}hmer and M.~Slassi.
\newblock Limit laws for distorted return time processes for infinite measure
  preserving transformations.
\newblock {\em arXiv:math. DS/0509609}, pages 1--20, 2005.

\bibitem[Sen76]{Seneta:76}
E.~Seneta.
\newblock {\em Regularly varying functions}.
\newblock Springer-Verlag, Berlin, 1976.
\newblock Lecture Notes in Mathematics, Vol. 508.

\bibitem[Tha95]{thaler:95}
M.~Thaler.
\newblock A limit theorem for the {P}erron-{F}robenius operator of
  transformations on {$[0,1]$} with indifferent fixed points.
\newblock {\em Israel J. Math.}, 91(1-3):111--127, 1995.

\bibitem[Tha98]{Thaler:98}
M.~Thaler.
\newblock The {D}ynkin-{L}amberti arc-sine laws for measure preserving
  transformation.
\newblock {\em Trans. Amer. Math. Soc.}, 350:4593--4607, 1998.

\bibitem[Tha00]{Thaler:00}
M.~Thaler.
\newblock The asymptotics of the {P}erron-{F}robenius operator of a class of
  interval maps preserving infinite measures.
\newblock {\em Studia Math.}, 143(2):103--119, 2000.

\bibitem[Tha05]{Thaler:05}
M.~Thaler.
\newblock Asymptotic distributions and large deviations for iterated maps with
  an indifferent fixed point.
\newblock {\em Stoch. Dyn.}, 5(3):425--440, 2005.

\bibitem[TZ06]{ThalerZweimueller:06}
M.~Thaler and R.~Zweim{\"u}ller.
\newblock Distributional limit theorems in infinite ergodic theory.
\newblock {\em Probab. Theory Related Fields}, 135(1):15--52, 2006.

\bibitem[Zwe03]{Zweimueller:03}
R.~Zweim{\"u}ller.
\newblock Stable limits for probability preserving maps with indifferent fixed
  points.
\newblock {\em Stoch. Dyn.}, 3(1):83--99, 2003.

\end{thebibliography}
\end{document}